\documentclass[12pt,twoside,a4paper,titlepage]{amsart}
\usepackage[ansinew]{inputenc}
\usepackage{comment}
\usepackage{amsmath}
\usepackage{amsfonts}
\usepackage{amssymb}
\usepackage{mathtools}
\usepackage[arrow, matrix, curve]{xy}
\usepackage{graphicx}
\usepackage[left=3cm,top=3cm,bottom=3.5cm]{geometry}
\usepackage{fancyhdr}
\usepackage{amsthm}
\usepackage{paralist}
\usepackage{tikz}
\usepackage{pgf}
\usepackage{pgfplots}
\usepackage{textcmds}
\usepackage{microtype}
\usepackage{hyperref}
\usepackage{xy}
\usetikzlibrary{shapes.geometric}
\newcommand{\D}[1]{\textit{\textbf{#1}}}
\usepackage{cleveref}
\usepackage{xcolor}
\usepackage{float}
\colorlet{mygreen}{black!50!green}
\theoremstyle{definition}
\newtheorem{Def}{Definition}[section]
\newtheorem{Ex}[Def]{Example}
\newtheorem{Not}[Def]{Notation}
\newtheorem{Con}[Def]{Construction}
\theoremstyle{plain}
\newtheorem{Lem}[Def]{Lemma}
\newtheorem{Theo}[Def]{Theorem}
\newtheorem{Prop}[Def]{Proposition}
\newtheorem{Cor}[Def]{Corollary}
\theoremstyle{remark}
\newtheorem{Rem}[Def]{Remark}
\newtheoremstyle{beweis}
{3pt}%					<freier Platz darÃƒÂ¼ber>
{3pt}%					<freier Platz darunter>
{default}%			<Schriftart Text>
{\parindent}%		<Einschub>
{default}%			<Schriftart Titel>
{:}%						<Zeichen zwischen Titel und Text>
{\newline}%			<Platz nach dem Titel>
{}%							<Titel Spezifikationen>
\theoremstyle{beweis}

\newtheoremstyle{TheoremNum}
{\topsep}{\topsep}              %%% space between body and thm
{\itshape}                      %%% Thm body font
{}                              %%% Indent amount (empty = no indent)
{\bfseries}                     %%% Thm head font
{.}                             %%% Punctuation after thm head
{ }                             %%% Space after thm head
{\thmname{#1}\thmnote{ \bfseries #3}}%%% Thm head spec
\theoremstyle{TheoremNum}
\newtheorem{thmn}{Theorem}
\crefname{Ex}{Example}{Examples}
\crefname{Prop}{Proposition}{Propositions}
\crefname{Rem}{Remark}{Remarks}
\crefname{Cor}{Corollary}{Corollaries}
\crefname{Lem}{Lemma}{Lemmas}
\crefname{Theo}{Theorem}{Theorems}
\crefname{Def}{Definition}{Definitions}
\crefname{Not}{Notation}{Notations}

\title{On Merge Trees and Discrete Morse Functions on Paths and Trees}
\author{Julian Br\"uggemann}
\address[Julian Br\"uggemann]{Max Planck Institute for Mathematics, Bonn, Germany}
\email{brueggemann@mpim-bonn.mpg.de}
\usetikzlibrary{patterns}
\begin{document}
    \begin{abstract}
	In this work we answer an open question asked by Johnson--Scoville. We show that each merge tree is represented by a discrete Morse function on a path. Furthermore, we present explicit constructions for two different but related kinds of discrete Morse functions on paths that induce any given merge tree. A refinement of the used methods allows us to define notions of equivalence of discrete Morse functions on trees which give rise to a bijection between equivalence classes of discrete Morse functions and isomorphism classes of certain labeled merge trees. We also compare our results to similar ones from the literature, in particular to work by Curry.
	\end{abstract}
	\begingroup
	\let\newpage\relax
	\maketitle
    \endgroup 
    \tableofcontents
	\section{Introduction}
	Discrete Morse theory is a combinatorial version of the classical smooth Morse theory. It was originally developed by Forman in \cite{Mor}.\par
    In discrete Morse theory, topological properties of simplicial complexes $X$ are analyzed by considering discrete Morse functions $f\colon X \rightarrow \mathbb{R}$. These topological properties can in turn be used to obtain cell decompositions of $X$ with fewer cells. A good introduction to the topic is found in \cite{User}.\par
	Merge trees are used in Morse theory in order to keep track of the development of connected components of sublevel sets %\footnote{Some authors use superlevel sets instead of sublevel sets which gives a different notion of persistent connectivity.}
	$X_a\coloneqq f^{-1}(-\infty,a]$ of a given Morse function $f\colon X\rightarrow \mathbb{R}$. Since the sublevel sets form a filtration of $X$, merge trees can be seen as a combinatorial description of the persistent connectivity of $X$. In particular, every branching in the induced merge tree $M(X,f)$ corresponds to a pair of connected components of a sublevel set $X_{a-\varepsilon}$ that merge to one connected component in a sublevel complex  $X_a$ of higher level. \par
	Initially, merge trees were introduced to topological data analysis as an approximation to the Reeb graph, respectively contour tree, in \cite{Contour}. The Reeb graph is a graph that keeps track of the connected components of level sets of any given filtered manifold. In applications, the data set is often interpreted as a sampling of the graph of a function rather than a more general manifold, which is why the Reeb graph is often actually a tree, the so-called contour tree. \par
	Among other implementations of techniques of smooth Morse theory, computational methods for the Reeb graph have originally been introduced in \cite{SurfCode} in order to handle surfaces embedded in 3D with the help of computers. Later on, several ways to compute and apply contour trees of data sets have been discussed in many articles, e.g. \cite{ExTo}, \cite{ConTreSma}, \cite{ConofCon}, and \cite{Contour}. Among other applications, merge trees have been used in visualization, e.g. in \cite{Landsc}, \cite{Landscmeta}, \cite{CaVTVMT}, and \cite{StrucAv}. Surveys about applications of merge trees and other concepts in visualization can be found in \cite{Survey1} and \cite{Survey2}. Furthermore, a certain version of chiral merge trees has been used in \cite{TSPHC} to analyze asymmetries of time series.\par  
	We focus on the more structural and theoretical side of merge trees, in particular the connection to discrete Morse theory. We consider a specific construction for merge trees induced by discrete Morse functions on trees which was introduced in \cite{JS}. We use this construction to gain a better understanding of the set of discrete Morse functions, the set of merge trees, and the relationship between the two. \par
	Similar work has been done in \cite{Fiber} for the relationship between Morse-like functions on the interval and a related version of merge trees. Furthermore, our work is in some sense similar to parts of \cite{TrBCaba}, with the difference being that the authors of \cite{TrBCaba} consider the relationship between merge trees and their induced barcodes instead. It seems reasonable to adapt the following terms from both of these articles and say that we consider a different instance of the ``fiber of the persistence map'', respectively a different instance of the ``inverse problem'' of the ``persistence map''. These names refer to the fact that merge trees and barcodes are invariants of filtered spaces rather than just spaces. In this context, the name ``persistent'' became popular due to persistent homology as it appears in topological data analysis. Similar to the authors of \cite{Fiber} and \cite{TrBCaba}, we are interested in finding out what information is lost by considering our invariant at hand, the induced merge tree, instead of the given data, in our case a discrete Morse function, and how this information might be re-obtained. Such knowledge might be helpful to investigate the space of merge trees and the space of discrete Morse functions in future work. Furthermore, a good understanding of the fiber of the persistence map is useful for topological data analysis because it hints at features which might be lost due to the chosen invariant. Moreover, insights about the inverse problem might be helpful to enhance the chosen invariant in a way such that it preserves certain desired features of the data set. \par
%	For instance, in \cite{JS} an equivalence relation between discrete Morse functions (dMf) on a given tree is introduced. Two dMfs are considered to be equivalent if and only if their induced merge trees are isomorphic. \par
	We respond to an open question asked in \cite{JS} by showing that every merge tree is represented by a discrete Morse function (dMf).
	In particular, for any given merge tree we construct a dMf on a path as a representative of the isomorphism class defined by said merge tree:\par
	\mbox{}\par
	\begin{thmn}[\ref{MainT}]
		Let $T$ be a merge tree. Then there is a path $P$ such that $T\cong M(P,f_{io})\cong M(P,f_{sc})$ holds as merge trees where $f_{io}$ denotes the induced index-ordered dMf (\Cref{f'}) and $f_{sc}$ denotes the sublevel-connected dMf (\Cref{subcondmf}) on $P$.
	\end{thmn}
	\mbox{}\par
	In particular, the discrete Morse function from \Cref{MainT} can be chosen to be index-ordered (\Cref{dmf}) or sublevel-connected (\Cref{defsubcon}). \par
	The main tool for the construction is the corresponding Morse order (\Cref{MO}), that is, the index Morse order (\Cref{MO}) or the sublevel-connected Morse order (\Cref{defsubcon}) on the nodes of a given merge tree $T$. The index Morse order defines leaf nodes to be strictly less than inner nodes. Among leaf nodes and among inner nodes, the index Morse order is defined by using a twisted version of length-lexicographical order on the set of path words (\Cref{Path}) that correspond to the respective nodes. The path words are defined by the chirality of the nodes of the shortest path from the root to the corresponding node. For the sublevel-connected Morse order (\Cref{scmo}) we do not artificially distinguish between leaf nodes and inner nodes.\par
    We use the index Morse order (\Cref{MO}) to define the index Morse labeling (\Cref{iML}) on the nodes of $T$. Together with the simplex order (\Cref{SO}), which establishes a correspondence (\Cref{iso}) between the nodes of $T$ and the simplices of a path $P$, the index Morse labeling defines the induced index-ordered discrete Morse function on said path $P$. \par
    In \Cref{slc} we introduce several kinds of equivalence relations on the sets of discrete Morse functions with only critical cells on paths and trees. These equivalence relations allow us to identify equivalence classes of discrete Morse functions with isomorphism classes of Morse labeled merge trees:
    \begin{thmn}[\ref{dmfml}]
    	The induced labeled merge tree $M(\_ \ , \_)$ and the induced dMf $\Phi$ define maps $M(\_ \ , \_)\colon DMF^{\text{crit}}_P \xleftrightarrow{} MlT \colon \Phi$ that are inverse to each other in the sense that:
    	\begin{enumerate}
    		\item for a dMf $(P,f)$ with only critical cells, the dMf $\Phi (M(P,f),\lambda_f)$ is symmetry-equivalent to $(P,f)$, and
    		\item for an Ml tree $(T,\lambda)$, the Ml tree $M(\Phi T,f_\lambda)$ is isomorphic to $(T,\lambda)$.
    	\end{enumerate}
    \end{thmn}

\begin{thmn}[\ref{dmfmlt}]
	The induced labeled merge tree $M(\_ \ , \_)$ and the induced dMf $\Phi$ define maps $M(\_ \ , \_)\colon DMF^{\text{crit}}_X \xleftrightarrow{} MlT \colon \Phi$ that are inverse to each other in the sense that:
	\begin{enumerate}
		\item for any dMf $(X,f)$ with only critical cells, the dMf $\Phi (M(X,f),\lambda_f)$ is cm-equivalent to $(X,f)$, and
		\item for any Ml tree $(T,\lambda)$, the Ml tree $M(\Phi T,f_\lambda)$ is isomorphic to $(T,\lambda)$.
	\end{enumerate}
\end{thmn}
    The construction of the discrete Morse function induced by a Morse labeling is similar to the construction of functions on the interval in \cite{Fiber}. In particular, \Cref{dmfml} is very similar to the result \cite[Prop 6.11]{Fiber}. The use of Morse labelings in this work basically plays the role of the function $\pi\colon T \rightarrow \mathbb{R}$ from \cite{Fiber}. Moreover, the simplex order is almost the same as the use of chirality in \cite[Lem 6.4]{Fiber}. But in this work, Morse orders, and in turn Morse labelings, have to satisfy a certain compatibility with the chirality, that is, property (2) of \Cref{morseorder}.\par 
    The notion of merge trees we use originates from \cite[Def 5]{JS} and differs from the one used in \cite{Fiber}: A priori, merge trees $T$ in the sense of \cite{JS} do not carry a height function $T \rightarrow \mathbb{R}$ as part of their data. Instead, the two children of each node have a chirality assigned to them as part of the tree's data. This means that for any two child nodes of the same parent node, it is specified as part of data which is the right and which is the left child. This version of chirality is also canonically assigned to merge trees induced by discrete Morse functions. In contrast, the chirality of chiral merge trees in the sense of \cite{Fiber} arises from a chosen orientation on the interval. We obtain a similar correspondence between the chirality of merge trees and orientations on paths using the simplex order, \Cref{SO}. Apart from these differences, the notion of merge trees in the sense of \cite{JS} is closely related to the one from \cite{Fiber}. Chirality in the sense of \cite{JS} is a specific version of the notion of chirality used in \cite{Fiber}. In order to see this, we show that the construction of the merge tree $M(X,f)$ induced by a discrete Morse function $f$ in the sense of \cite{JS} can be modified, see \Cref{indMltree}, to obtain a function $T \rightarrow \mathbb{R}$ from $f$, similarly to \cite{Fiber}. This gives rise to the notion of Morse labelings, \Cref{ML}, and Morse orders, \Cref{MO}. It turns out that the use of chirality in \cite{JS} assumes a certain compatibility between Morse orders and the simplex order, whereas the use of chirality in \cite{Fiber} does not. As a result, the induced merge tree in the sense of \cite{Fiber} distinguishes between symmetry equivalences, \Cref{symm}, whereas the induced merge tree in the sense of \cite{JS} identifies symmetry-equivalent discrete Morse functions with each other, see \Cref{symequiv}. We discuss this in a bit more detail at the end of \Cref{Rel}. \par
    In said discussion, we mention the notion of CMl trees, see \Cref{CML}, which is as objects basically the same as the notion of merge trees from \cite[Def 2.2]{TrBCaba}. However, the notion of combinatorial equivalence of labeled merge trees from \cite[Def 2.6]{TrBCaba} corresponds to a non-chiral version of shuffle equivalence \Cref{MlTiso} rather than an equivalence of the induced persistent set as in \cite[5.1]{Fiber}, which is more similar to an isomorphism of Ml trees \Cref{MoTiso}. The chiral merge trees from \cite[Def 2.1]{TSPHC} are as objects also very similar to the chiral merge trees from \cite{Fiber} and, thus, differ similarly from our notion of labeled merge trees. \par
    The aforementioned versions of labeled merge trees all have in common that their labelings need to be compatible with some other data inherent to the merge tree. In contrast to that, the labelings from \cite{StrucAv} can be quite arbitrary and might even assign multiple labels to a single node. Hence, our Ml trees a priori seem to be unrelated to the notion of labeled merge trees from \cite{StrucAv}. \par
%	The step-by-step approach (\Cref{Con}) gives an alternative construction of the path $P$ and the induced index-ordered discrete Morse function on $P$ by induction over the nodes of $T$ from highest to lowest with respect to the index Morse order.\par
%	\Cref{equal} shows that the two constructions of the induced index-ordered discrete Morse functions yield the same result.\par
%	We illustrate the step-by-step construction in an example in \Cref{Ex}.\par
	
	\addtocontents{toc}{\protect\setcounter{tocdepth}{1}}
	\subsection*{Acknowledgements}
	The author would like to thank Benjamin Johnson and Nicholas A. Scoville for the suggested questions in their article \cite{JS} which led to this project. Furthermore, the author would like to thank the mathematical faculty of Ruhr-University of Bochum, especially the chair for topology, for the great scientific environment in which this endeavor was started. In addition to that the author thanks Max Planck Institute for Mathematics for the great scientific environment in which this project was finished.\par
	Most notably the author thanks his advisor, Viktoriya Ozornova, for her advice and the many helpful discussions which increased the quality of this article. Last but not least, the author thanks the anonymous referees for the helpful and detailed feedback.
	\section{Conflict of Interest Statement}
	The author states that there is no conflict of interest.
	\section{Preliminaries}
	\addtocontents{toc}{\protect\setcounter{tocdepth}{3}}
	We consider discrete Morse functions (dMf) on trees. Recall that trees are finite acyclic simple graphs. Furthermore, simple graphs are 1-dimensional simplicial complexes. Where feasible, we introduce the preliminaries in the broader generality they are usually defined in, rather than in the lesser generality we actually need for this work. We adapt most notations and conventions from \cite{JS}. For simplicity, we assume all trees in this article to be non-empty. 
%	Even though the empty cases are well-defined, they are not very interesting.
	Similar to \cite{JS}, we assume the dMfs to fulfill certain generic properties. In detail this means the following :\par
	\begin{Def}\label{dmf}
		Let $X$ be a simplicial complex. A map $f \colon X \rightarrow \mathbb{R}$ is called a \D{discrete Morse function (dMf)} if it fulfills the following properties for any pair of simplices $\sigma,\tau\in X$:
		\begin{enumerate}[(i)]
			\item $\sigma \subseteq \tau \Rightarrow f(\sigma) \leq f(\tau)$ (weakly increasing)
			\item  $f$ is at most $2-1$ 
			\item $f(\sigma)=f(\tau) \Rightarrow (\sigma \subset \tau \vee \tau \subset \sigma)$ (matching)
		\end{enumerate}
	Simplices on which $f$ is $1-1$ are called \D{critical}. Simplices which belong to the preimage of the same value are called \D{matched}. The set of critical simplices is denoted by $\operatorname{Cr}(f)$. Values of critical simplices under $f$ are called \D{critical values} of $f$. A dMf is called \D{index-ordered} if for arbitrary critical simplices $\sigma,\tau$ the following holds: If $\operatorname{dim}(\sigma)$ is smaller than $\dim(\tau)$, then $f(\sigma) < f(\tau)$ holds.
	\end{Def}
\begin{Rem}
	The definition given above is not the most general definition of dMfs but rather assumes several generic properties. This means that any dMf in the sense of \cite{Mor} can be modified by a Forman equivalence (see \cite[Def 4.1]{JS}) to fulfill these properties, that is, without changing the induced Morse matching. As usual in the context of dMfs, we write $f\colon X \rightarrow \mathbb{R}$ although the map $f$ is actually defined on the face poset of $X$. \par
	In \cite{NTT}, a similar generic property is used to analyze flow paths induced by dMfs. The notion of faithful dMfs as defined in \cite[Def 2.9]{NTT} is almost the same as index-ordered dMfs in this article. The only difference is that for index-ordered dMfs matched cells have the same value, whereas for faithful dMfs the values of the matched boundary simplices are higher than the values on the corresponding matched co-boundary simplices.\par
	Furthermore, since simple graphs and in particular trees are examples of simplicial complexes, the definition of dMfs can be applied to them as well. \par
\end{Rem}
\begin{Rem}\label{critorder}
    Because $f$ maps $\operatorname{Cr}(f)$ injectively into a totally ordered set, it induces a total order on $\operatorname{Cr}(f)$. We refer to this induced order whenever we speak of simplices being \D{ordered by $f$}.
\end{Rem}
\begin{Not}
		We use the following conventions regarding notation:
	\begin{itemize}
		\item We depict dMfs on graphs by labeling the graph with the values of the dMf. 
		\item Let $G$ be a graph and let $v$ be a node of $G$. By $G[v]$ we denote the connected component of $G$ which contains $v$.
	\end{itemize}
\end{Not}
\begin{Def}
\begin{itemize}
		\item Let $X$ be a simplicial complex, $f \colon X\rightarrow \mathbb{R}$ a dMf and $a\in \mathbb{R}$. The \D{sublevel complex} of level $a$, denoted by $X^f_a$, is defined by $X^f_a\coloneqq \{\sigma \in X \ \rvert \ f(\sigma) \leq a \}$. If the referred dMf $f$ is clear from the context, we drop the superscript $f$ from the notation.
%		the subcomplex of $X$ that consists of all simplices $\sigma$ with $f(\sigma)\leq a$.
		\item The ordered critical values $c_0< c_1< \dots < c_{m}$ induce a chain of sublevel complexes $X^f_{c_0}\subsetneq X^f_{c_1}\subsetneq \dots \subsetneq X^f_{c_{m}}$. Within this chain, we refer by $X^f_{c_i-\varepsilon}$ to the complex that immediately precedes $X_{c_i}$.
	\end{itemize}
\end{Def}
\begin{Rem}
	The given definition of sublevel complexes differs from the standard one used in the literature. We make use of the fact that a dMf $f$ being weakly increasing implies that $X_a^f$ as defined above is already a subcomplex of $X$. If we wanted to consider the general definition of dMfs as introduced in \cite{Mor}, we would have to work with the smallest supercomplex of $X_a^f$ in $X$ instead. Taking the smallest supercomplex corresponds to additionally including all faces of simplices of $X_a^f$ to make $X_a^f$ a simplicial complex.  
%	The aforementioned notion of sublevel complex is only that concise because dMfs are already assumed to be weakly increasing in this paper. If you work with a more general notion of dMf, you will need to include all faces of simplices $\sigma$ with $f(\sigma) \leq a$ as well.
\end{Rem}
\begin{Lem}\label{min}
    Let $X$ be a finite simplicial complex and let $f \colon X \rightarrow \mathbb{R}$ be a dMf. Then $f$ attains its minimum on a critical 0-simplex. Furthermore, the statement also holds for the restriction to any connected component of sublevel complexes.
\end{Lem}
\begin{proof}[Sketch of Proof]The statement follows by a proof by contradiction and properties (i) and (ii) of \Cref{dmf}.
%    Because $f$ is weakly increasing, no higher simplex $\tau$ can have a value which is strictly smaller than the value of $f$ on any of $\tau$s boundary 0-simplices. Thus, the minimum is attained on a 0-simplex. \par
%    Let $\sigma$ be a 0-simplex such that $f(\sigma$) is the minimum of $f$. We prove by contradiction that $\sigma$ is critical: If $\sigma$ is not critical, then there is a 1-simplex $\tau$ such that $\sigma$ is a boundary 0-simplex of $\tau$ and $f(\tau)=f(\sigma)$ holds. Let $\sigma'$ be the other boundary 0-simplex of $\tau$. Then $f(\sigma')<f(\tau)=f(\sigma)$ holds because $f$ is by assumption weakly increasing and at most $2-1$. This is a contradiction to $f$ attaining its minimum on $\sigma$, so $f$ attains its minimum on a critical 0-simplex.\par
%    Since all three properties of \Cref{dmf} are inherited by arbitrary subcomplexes of $X$, the statement is also true for connected components of sublevel complexes.
\end{proof}
\begin{Rem}
    The analogous statement for the maximum of a dMf $f$ on arbitrary 1-simplices is false, as the following example shows:\\
    \begin{tikzpicture}
        \draw (-6.5,1.5) node {};
        \draw (0,1) node {$\bullet$};
        \draw (1,1) node {$\bullet$};
        \draw (2,1) node {$\bullet$};
        \draw (0,1) -- (2,1);
        \draw (0,1.3) node {$3$};
        \draw (1,1.3) node {$0$};
        \draw (2,1.3) node {$1$};
        \draw (0.5,1.3) node {$3$};
        \draw (1.5,1.3) node {$2$};
    \end{tikzpicture}\par
    Here, the maximum is attained on a pair of matched simplices.
\end{Rem}
	\subsection{Merge Trees}\label{subsecmerge}
	We briefly recapture preliminaries about merge trees as they are explained in \cite{JS}. The basic idea is that merge trees keep track of the chronological development of the connected components of sublevel sets. We adapt the point of view of \cite{JS}, that is, we consider them \q{upside down}. Thus, the children will appear above their parent node. Afterwards we introduce additional structure that dMfs induce on their corresponding merge trees and consider notions of equivalence which arise from that structure. 
	\begin{Def}[Merge Tree]
		A \D{merge tree} is a full rooted chiral binary tree $T$. In detail this means that $T$ is a rooted tree fulfilling the properties of being binary and full, and that $T$ has the extra datum of being chiral: 
		\begin{enumerate}
			\item[\textbf{full binary:}] Each node of $T$ has either zero or two children.
			\item[\textbf{chiral:}] Each child node in $T$ carries the extra datum, the so-called \D{chirality}, of being a left or a right child.
		\end{enumerate}
		  \D{Morphisms} of merge trees are morphisms of rooted binary trees which are compatible with the chirality.\par
		\end{Def}
		For rooted trees $T$ we use the notions of \D{subtrees}, \D{ancestors} and \D{descendants} as they are commonly used in computer science. 
		\begin{Def}\label{subtree}
		For any node $p$ of $T$, the \D{descendants} of $p$ are defined inductively: A node $c$ is a descendant of $p$ if the parent node of $c$ is a descendant of $p$ or $p$ itself.\par 
		A \D{subtree} of $T$ is a subgraph of $T$ that consists of exactly all of the descendants of some node $p$ of $T$. \par
		For a node $p$ of $T$ we call all nodes which lie on the shortest path between $p$ and the root, including the root, the \D{ancestors} of $p$.
		\end{Def}
		\begin{Not}\label{Nota}
	For an inner node $c$ of $T$, we denote the left child of $c$ with $c_l$ and the right child of $c$ with $c_r$. We illustrate this notation in the following example:\\
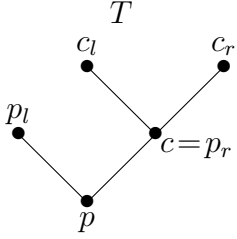
	\begin{tikzpicture}[scale=0.9]
	\draw (-6,-1) node {};
	\draw (2.5,-1.2) node {$T$};
	\draw (2,-4) node {$\bullet$};
	\draw (2,-4.3) node {$p$};
	\draw (1,-3) node {$\bullet$};
	\draw (1,-2.7) node {$p_l$};
	\draw (3,-3) node {$\bullet$};
	\draw (3.6,-3.2) node {$c \! =\! p_r$};
	\draw (2,-2) node {$\bullet$};
	\draw (2,-1.7) node {$c_l$};
	\draw (4,-2) node {$\bullet$};
	\draw (4,-1.7) node {$c_r$};
	\draw (2,-4) -- (4,-2);
	\draw (2,-4) -- (1,-3);
	\draw (3,-3) -- (2,-2);
	\draw (2,-4) node {};
	\end{tikzpicture}
	\end{Not}
		\begin{Rem}\label{BinTree}
	For full binary trees $T$ with $i(T)$ inner nodes and $l(T)$ leaves it is a well-known result that $l(T)=i(T)+1$ holds. It can be proved inductively. 
\end{Rem}
		\begin{Rem}
			The chirality of nodes will either be denoted by labels or indicated implicitly by embedding the merge tree on the page. Throughout the literature there are different notions of merge trees that are not always distinguished by name or notation. In this work, merge trees do not have explicit weights on edges. Moreover, merge trees in this work a priori do not carry a function to the real numbers. In that way, we distinguish between merge trees and Morse labeled merge trees which will be introduced in \Cref{ML}. We use the chirality to obtain Morse labelings on unlabeled merge trees. This leads to \Cref{MO}.
		\end{Rem}
		
		\begin{Con}[{\cite[Thm 9]{JS}}]\label{Mer}
		Let $X$ be a tree and let $f\colon X \rightarrow \mathbb{R}$ be a dMf. The \D{merge tree induced by $f$}, denoted by $M(X,f)$ is constructed as follows:\par
		Let $c_0<c_1<\dots < c_m$ be the critical values of $f$ that are assigned to 1-simplices. The associated merge tree $M(X,f)$ is constructed by induction over these critical values in descending order. Furthermore, we label the nodes of $M(X,f)$ in order to refer to them later. The label of a node $n$ will be denoted by $\lambda(n)$.\par
		For the base case we begin by creating a node $M(c_m)$ which corresponds to the critical 1-simplex in $X$ labeled $c_m$ and setting its label $\lambda(M(c_m))$ to $(c_m ,L)$. \par
		For the inductive step, let $M(c_i)$ be a node of $M(X,f)$ that corresponds to a critical 1-simplex between two 0-simplices $v$ and $w$. Define $\lambda_v\coloneqq$ $\max\{f(\sigma) \rvert \sigma \in X_{c_i-\varepsilon}[v], \sigma \text{ critical}\}$ and $\lambda_w\coloneqq\max\{f(\sigma) \vert \sigma \in X_{c_i-\varepsilon}[w], \sigma \text{ critical}\}$.
		Two child nodes of $M(c_i)$ are created, named $n_{\lambda_v}$ and $n_{\lambda_w}$. 
		Then label the new nodes $\lambda(n_{\lambda_v})\coloneqq\lambda_v$ and $\lambda(n_{\lambda_w})\coloneqq\lambda_w$. If $\min\{f(\sigma) \vert \sigma \in X_{c_i-\varepsilon}[v]\}<\min\{f(\sigma) \vert \sigma \in X_{c_i-\varepsilon}[w]\}$, we assign $n_{\lambda_v}$ the same chirality (L or R) as $M(c_i)$ and give $n_{\lambda_w}$ the opposite chirality. Continue the induction over the rest of the critical 1-simplices.
	\end{Con}
\begin{Rem}
	By construction, one of the following two cases holds for the labels $\lambda_v$ and $\lambda_w$.
	They might be either critical values lower than $c_i$ that are assigned to edges or critical values that are assigned to nodes. The two labels $\lambda_v$ and $\lambda_w$ do not necessarily belong to the same case.\par
	Therefore, the nodes $n_{\lambda_v}$ and $n_{\lambda_w}$ will possibly be denoted as $M(c_j)$ and $M(c_k)$ for some $j,k<i$ in later steps of the induction. In particular, this means that the node which is considered in the first instance of the inductive step is $c_m$.
\end{Rem}
\begin{Rem}
    Although the induced merge tree as introduced above comes with a labeling, the labeling is not part of the data of the induced merge tree. This is one of the main differences between merge trees in \cite{JS} and the merge trees in \cite{Fiber}.
\end{Rem}
We consider one way to keep some information provided by the induced labeling on $M(X,f)$.
\begin{Def}\label{morseorder}
     Let $T$ be a merge tree. We call a total order $\leq$ on the nodes of $T$ a \D{Morse order} if it fulfills the following two properties for any subtree $T'$ of $T$:
     \begin{enumerate}
        \item The restriction $\leq_{\lvert T'}$ attains its maximum on the root of $T'$.
        \item The restriction $\leq_{\lvert T'}$ attains its minimum on the subtree with root $\textcolor{red}{p_{l}}/\textcolor{blue}{p_r}$ if \textcolor{red}{L}/\textcolor{blue}{R} is the chirality of the root $p$ of $T'$.
    \end{enumerate}
    Moreover, we call a merge tree $(T,\leq)$ together with a Morse order $\leq$ a \D{Morse-ordered merge tree (Mo tree)}.
\end{Def}
\begin{Rem}
	Assuming property (2) of \Cref{morseorder} for every subtree $T'$ with root $p$ of $T$ is equivalent to either of the following:
	\begin{itemize}
		\item For any subtree $T'$ with root $p$ of $T$, the minimum of $\leq_{\lvert T'}$ has the same chirality as $p$.
		\item For any subtree $T'$ with root $p$ of $T$, all nodes on the shortest path between $p$ and the minimum of $\leq_{\lvert T'}$ have the same chirality as $p$. 
	\end{itemize}
	The equivalence can be proved by an inductive argument over all nodes of the shortest path between $p$ and the minimum.
\end{Rem}
Because Morse orders $\leq$ define in particular finite totally ordered sets, there are unique order-preserving isomorphisms \mbox{$\lambda \colon (V(T),\leq) \xrightarrow[]{\cong} \{0,1,\dots, i(T)+l(T)-1\} \subseteq \mathbb{N}_0\subset \mathbb{R}$} for each Morse order. Conversely, each injective labeling $\lambda \colon T \rightarrow \mathbb{R}$ induces an order $\leq_\lambda$ on the nodes of $T$ by usage of the total order on $\mathbb{R}$.
\begin{Def}\label{ML}
     For Morse orders $\leq$ we call the map $\lambda_{\leq} \colon (V(T),\leq) \rightarrow \{0,1,\dots, i(T)+l(T)-1\}$ the \D{Morse labeling induced by $\leq$}. We call an arbitrary labeling $\lambda \colon T \rightarrow \mathbb{R}$ a \D{Morse labeling} if it induces a Morse order on $T$.\par
     We call a merge tree $(T,\lambda)$ with a Morse labeling $\lambda \colon T \rightarrow \mathbb{R}$ a \D{Morse labeled merge tree (Ml tree)}.\par
     For a Mo tree $(T,\leq)$ we call the Ml tree $(T,\lambda_{\leq})$ the \D{Ml tree induced by $(T,\leq)$}.
\end{Def}
%\begin{Rem}
%\end{Rem}
\begin{Prop}\label{indMltree}
    Let $f\colon X \rightarrow \mathbb{R}$ be a dMf. The labeling which appears in \Cref{Mer} induces a Morse order on $M(X,f)$. Hence, $M(X,f)$ canonically carries the structure of a Mo tree as well as an Ml tree.
\end{Prop}
\begin{proof}
    It is proved in \cite[Thm 9]{JS} that $M(X,f)$ is a merge tree. We only have to prove that the labeling induces a Morse order. By \Cref{critorder} the set $\operatorname{Cr}(f)$ of critical values carries a total order induced by $f$. Since the critical values of $f$ precisely define the labeling in \Cref{Mer}, the labeling induces a total order on the nodes of $M(X,f)$. It is only left to prove that this order is a Morse order.\par
    In the construction, each inner node of $M(X,f)$ corresponds to a critical 1-simplex and is labeled with the critical value of said critical 1-simplex. Since parent nodes are created before their child nodes are, and since the critical values are considered from highest to lowest, property (1) of \Cref{morseorder} is fulfilled. The rule in the construction which decides the chirality of the child nodes is exactly the same as property (2) of \Cref{morseorder}. Hence, it is fulfilled by construction.
\end{proof}
We denote the canonical labeling of $M(X,f)$ by $\lambda_f$.
%\begin{Rem}\label{commute}
%It is immediate that the underlying merge tree of the Ml tree induced by a dMf $f$ is just the merge tree induced by $f$ in the sense of \Cref{Mer}. That is, taking the induced merge tree commutes with taking the induced Ml tree and forgetting the Morse label.
%\end{Rem}
\begin{Rem}
    In the aforementioned proof, it becomes clear that both conditions of \Cref{morseorder} are necessary for a total order on $M(X,f)$ to be induced by a dMf $f$. \par
    Morse orders are useful for the construction of dMfs that induce given merge trees $T$.
    We will see in \Cref{induceddmfisdmf} that for a total order on an arbitrary merge tree $T$ condition (1) is sufficient for inducing a dMf in the sense of \Cref{f'} later on. Nonetheless, condition (2) is necessary to ensure that the induced dMf induces the given merge tree $T$.  
\end{Rem}
\begin{Ex}
	We consider the following example of a dMf $f \colon X\rightarrow \mathbb{R}$:
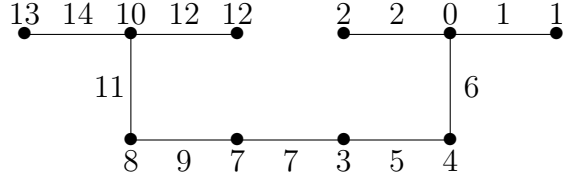
\par
	\begin{tikzpicture}[scale=0.7]
	\draw (-3-1.9,3) node {};
	\draw (0,2) node {$\bullet$};
	\draw (2,2) node {$\bullet$};
	\draw (4,2) node {$\bullet$};
	\draw (2,0) node {$\bullet$};
	\draw (4,0) node {$\bullet$};
	\draw (6,0) node {$\bullet$};
	\draw (8,0) node {$\bullet$};
	\draw (6,2) node {$\bullet$};
	\draw (8,2) node {$\bullet$};
	\draw (10,2) node {$\bullet$};
	\draw (0,2) -- (4,2);
	\draw (6,2) -- (10,2);
	\draw (2,0) -- (2,2);
	\draw (2,0) -- (8,0);
	\draw (8,0) -- (8,2);
	\draw (0,2+0.4) node {13};
	\draw (2,2+0.4) node {10};
	\draw (4,2+0.4) node {12};
	\draw (2,0-0.4) node {8};
	\draw (4,0-0.4) node {7};
	\draw (6,0-0.4) node {3};
	\draw (8,0-0.4) node {4};
	\draw (6,2+0.4) node {2};
	\draw (8,2+0.4) node {0};
	\draw (10,2+0.4) node {1};
	\draw (9,2+0.4) node {1};
	\draw (7,2+0.4) node {2};
	\draw (3,2+0.4) node {12};
	\draw (1,2+0.4) node {14};
	\draw (7,0-0.4) node {5};
	\draw (5,0-0.4) node {7};
	\draw (3,0-0.4) node {9};
	\draw (1.6,1) node {11};
	\draw (8.4,1) node {6};
	\end{tikzpicture}\\
	The critical values on edges are:
	\[
	5<6<9<11<14
	\]
	We now show the construction algorithm of $M(X,f)$ visually by depicting $X_{c_i-\varepsilon}$ on the left and the part of $M(X,f)$ that is created up to the step corresponding to $c_i$ on the right.\\
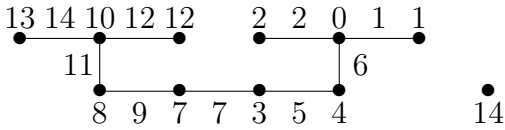
	\begin{tikzpicture}[scale=0.7]
%	\draw (0-1.2,100) node {};
	\draw (1,99.3) node {Start:};
	\draw (1,98) node {$\bullet$};
	\draw (2.5,98) node {$\bullet$};
	\draw (4,98) node {$\bullet$};
	\draw (5.5,98) node {$\bullet$};
	\draw (7,98) node {$\bullet$};
	\draw (8.5,98) node {$\bullet$};
	\draw (2.5,97) node {$\bullet$};
	\draw (4,97) node {$\bullet$};
	\draw (5.5,97) node {$\bullet$};
	\draw (7,97) node {$\bullet$};
	\draw (1,98) -- (4,98);
	\draw (5.5,98) -- (8.5,98);
	\draw (2.5,97) -- (7,97);
	\draw (2.5,98) -- (2.5,97);
	\draw (7,98) -- (7,97);
	\draw (1,98+0.4) node {13};
	\draw (1.75,98+0.4) node {14};
	\draw (3.25,98.4) node {12};
	\draw (2.5,98+0.4) node {10};
	\draw (4,98+0.4) node {12};
	\draw (5.5,98+0.4) node {2};
	\draw (6.25,98.4) node {2};
	\draw (7,98+0.4) node {0};
	\draw (7.75,98.4) node {1};
	\draw (8.5,98+0.4) node {1};
	\draw (2.5,97-0.4) node {8};
	\draw (3.25,96.6) node {9};
	\draw (4,97-0.4) node {7};
	\draw (4.75,96.6) node {7};
	\draw (5.5,97-0.4) node {3};
	\draw (6.25,96.6) node {5};
	\draw (7,97-0.4) node {4};
	\draw (2.1,97.5) node {11};
	\draw (7.4,97.5) node {6};
	\draw (9.8,97) node {$\bullet$};
	\draw (9.8,96.6) node {14};
	\end{tikzpicture}
	
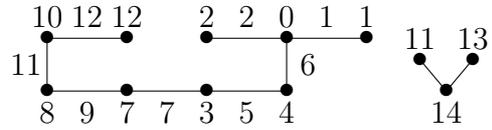
\begin{tikzpicture}[scale=0.7]
	\draw (0-0.2,96) node {};
	\draw (1,95.3) node {$c_i=14$:};
	\draw (1,98-4) node {$\bullet$};
	\draw (2.5,98-4) node {$\bullet$};
	\draw (4,98-4) node {$\bullet$};
	\draw (5.5,98-4) node {$\bullet$};
	\draw (7,98-4) node {$\bullet$};
	\draw (8.5,98-4) node {$\bullet$};
	\draw (2.5,97-4) node {$\bullet$};
	\draw (4,97-4) node {$\bullet$};
	\draw (5.5,97-4) node {$\bullet$};
	\draw (7,97-4) node {$\bullet$};
	\draw (2.5,98-4) -- (4,98-4);
	\draw (5.5,98-4) -- (8.5,98-4);
	\draw (2.5,97-4) -- (7,97-4);
	\draw (2.5,98-4) -- (2.5,97-4);
	\draw (7,98-4) -- (7,97-4);
	\draw (1,98+0.4-4) node {13};
	\draw (3.25,98.4-4) node {12};
	\draw (2.5,98+0.4-4) node {10};
	\draw (4,98+0.4-4) node {12};
	\draw (5.5,98+0.4-4) node {2};
	\draw (6.25,98.4-4) node {2};
	\draw (7,98+0.4-4) node {0};
	\draw (7.75,98.4-4) node {1};
	\draw (8.5,98+0.4-4) node {1};
	\draw (2.5,97-0.4-4) node {8};
	\draw (3.25,96.6-4) node {9};
	\draw (4,97-0.4-4) node {7};
	\draw (4.75,96.6-4) node {7};
	\draw (5.5,97-0.4-4) node {3};
	\draw (6.25,96.6-4) node {5};
	\draw (7,97-0.4-4) node {4};
	\draw (2.1,97.5-4) node {11};
	\draw (7.4,97.5-4) node {6};
	\draw (13-3,97-4) node {$\bullet$};
	\draw (13-3,96.6-4) node {14};
	\draw (12.5-3,96.6-3) node {$\bullet$};
	\draw (13.5-3,96.6-3) node {$\bullet$};
	\draw (12.5-3,97-3) node {11};
	\draw (13.5-3,97-3) node {13};
	\draw (13-3,97-4) -- (12.5-3,96.6-3);
	\draw (13-3,97-4) -- (13.5-3,96.6-3);
	\end{tikzpicture}\\
	
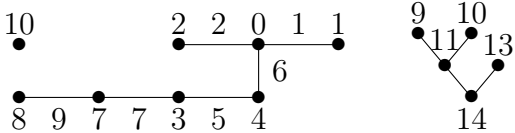
\begin{tikzpicture}[scale=0.7]
%	\draw (0-1.2,92) node {};
	\draw (2.5,91.3) node {$c_i=11$:};
	\draw (2.5,98-4-4) node {$\bullet$};
	\draw (5.5,98-4-4) node {$\bullet$};
	\draw (7,98-4-4) node {$\bullet$};
	\draw (8.5,98-4-4) node {$\bullet$};
	\draw (2.5,97-4-4) node {$\bullet$};
	\draw (4,97-4-4) node {$\bullet$};
	\draw (5.5,97-4-4) node {$\bullet$};
	\draw (7,97-4-4) node {$\bullet$};
	\draw (5.5,98-4-4) -- (8.5,98-4-4);
	\draw (2.5,97-4-4) -- (7,97-4-4);
	\draw (7,98-4-4) -- (7,97-4-4);
	\draw (2.5,98+0.4-4-4) node {10};
	\draw (5.5,98+0.4-4-4) node {2};
	\draw (6.25,98.4-4-4) node {2};
	\draw (7,98+0.4-4-4) node {0};
	\draw (7.75,98.4-4-4) node {1};
	\draw (8.5,98+0.4-4-4) node {1};
	\draw (2.5,97-0.4-4-4) node {8};
	\draw (3.25,96.6-4-4) node {9};
	\draw (4,97-0.4-4-4) node {7};
	\draw (4.75,96.6-4-4) node {7};
	\draw (5.5,97-0.4-4-4) node {3};
	\draw (6.25,96.6-4-4) node {5};
	\draw (7,97-0.4-4-4) node {4};
	\draw (7.4,97.5-4-4) node {6};
	\draw (13-2,89) node {$\bullet$};
	\draw (13-2,88.6) node {14};
	\draw (12.5-2,89.6) node {$\bullet$};
	\draw (13.5-2,96.6-3-4) node {$\bullet$};
	\draw (12-2,90.2) node {$\bullet$};
	\draw (12-2,90.6) node {9};
	\draw (13-2,90.2) node {$\bullet$};
	\draw (13-2,90.6) node {10};
	\draw (12.5-2,96.6-3-4) -- (13-2,90.2);
	\draw (12.5-2,97-3-3.9) node {11};
	\draw (13.5-2,97-3-4) node {13};
	\draw (13-2,97-4-4) -- (12-2,90.2);
	\draw (13-2,97-4-4) -- (13.5-2,96.6-3-4);
	\end{tikzpicture}
	
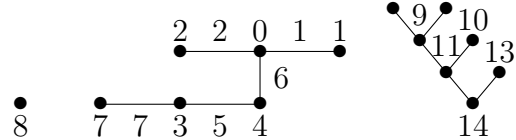
\begin{tikzpicture}[scale=0.7]
	\draw (1.2,88) node {};
	\draw (2,87.3) node {$c_i=9$:};
	\draw (5.5,98-12) node {$\bullet$};
	\draw (7,98-12) node {$\bullet$};
	\draw (8.5,98-12) node {$\bullet$};
	\draw (2.5,97-12) node {$\bullet$};
	\draw (4,97-12) node {$\bullet$};
	\draw (5.5,97-12) node {$\bullet$};
	\draw (7,97-12) node {$\bullet$};
	\draw (5.5,98-12) -- (8.5,98-12);
	\draw (4,97-12) -- (7,97-12);
	\draw (7,98-12) -- (7,97-12);
	\draw (5.5,86.4) node {2};
	\draw (6.25,98.4-12) node {2};
	\draw (7,98+0.4-12) node {0};
	\draw (7.75,98.4-12) node {1};
	\draw (8.5,98+0.4-12) node {1};
	\draw (2.5,97-0.4-12) node {8};
	\draw (4,97-0.4-12) node {7};
	\draw (4.75,96.6-12) node {7};
	\draw (5.5,97-0.4-12) node {3};
	\draw (6.25,96.6-12) node {5};
	\draw (7,97-0.4-12) node {4};
	\draw (7.4,97.5-12) node {6};
	\draw (13-2,89-4) node {$\bullet$};
	\draw (13-2,88.6-4) node {14};
	\draw (12.5-2,89.6-4) node {$\bullet$};
	\draw (13.5-2,96.6-3-4-4) node {$\bullet$};
	\draw (12-2,90.2-4) node {$\bullet$};
	\draw (12-2,90.7-4) node {9};
	\draw (13-2,90.2-4) node {$\bullet$};
	\draw (13-2,90.6-4) node {10};
	\draw (11.5-2,90.2-4+0.6) node {$\bullet$};
	\draw (11.5-2,86.2+1) node {6};
	\draw (12.5-2,85.8+1) node {$\bullet$};
	\draw (12.5-2,86.2+1) node {8};
	\draw (12.5-2,96.6-3-4-4) -- (13-2,90.2-4);
	\draw (12.5-2,97-3-4-3.9) node {11};
	\draw (13.5-2,97-3-4-4) node {13};
	\draw (13-2,97-12) -- (11.5-2,90.2-4+0.6);
	\draw (13-2,84+1) -- (13.5-2,84.6+1);
	\draw (12-2,85.2+1) -- (12.5-2,85.8+1);
	\end{tikzpicture}\\
	
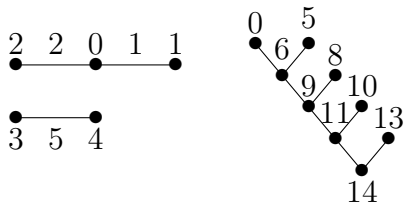
\begin{tikzpicture}[scale=0.7]
%	\draw (0-1.2,84) node {};
	\draw (4-2,82.3) node {$c_i=6$:};	
	\draw (5.5-2,98-16-1) node {$\bullet$};
	\draw (7-2,98-16-1) node {$\bullet$};
	\draw (8.5-2,98-16-1) node {$\bullet$};
	\draw (5.5-2,97-16-1) node {$\bullet$};
	\draw (7-2,97-16-1) node {$\bullet$};
	\draw (5.5-2,98-16-1) -- (8.5-2,98-16-1);
	\draw (5.5-2,97-16-1) -- (7-2,97-16-1);
	\draw (5.5-2,98+0.4-16-1) node {2};
	\draw (6.25-2,98.4-16-1) node {2};
	\draw (7-2,98+0.4-16-1) node {0};
	\draw (7.75-2,98.4-16-1) node {1};
	\draw (8.5-2,98+0.4-16-1) node {1};
	\draw (5.5-2,97-0.4-16-1) node {3};
	\draw (6.25-2,96.6-16-1) node {5};
	\draw (7-2,97-0.4-16-1) node {4};
	\draw (13-1-2,89-4-1-4-1) node {$\bullet$};
	\draw (13-1-2,88.6-4-1-4-1) node {14};
	\draw (12.5-1-2,89.6-4-1-4-1) node {$\bullet$};
	\draw (13.5-1-2,96.6-3-4-4-1-4-1) node {$\bullet$};
	\draw (12-1-2,90.2-4-1-4-1) node {$\bullet$};
	\draw (12-1-2,90.6-4-1-4-1) node {9};
	\draw (13-1-2,90.2-4-1-4-1) node {$\bullet$};
	\draw (13-1-2,90.6-4-1-4-1) node {10};
	\draw (11.5-1-2,90.2-4-1+0.6-4-1) node {$\bullet$};
	\draw (11.5-1-2,86.3-4-1) node {6};
	\draw (12.5-1-2,85.8-4-1) node {$\bullet$};
	\draw (12.5-1-2,86.2-4-1) node {8};
	\draw (11-1-2,90.2-4-1+0.6-4-1+0.6) node {$\bullet$};
	\draw (11-1-2,90.2-4-1+0.6-4-1+1) node {0};
	\draw (12.5-1-2,96.6-3-4-4-1-4-1) -- (13-1-2,90.2-4-1-4-1);
	\draw (12.5-1-2,97.1-3-4-4-1-4-1) node {11};
	\draw (13.5-1-2,97-3-4-4-1-4-1) node {13};
	\draw (11-2,90.2-4-1+0.6-4-1+0.6) node {$\bullet$};
	\draw (11-2,90.3-4-1+0.6-4-1+1) node {5};
	\draw (13-1-2,97-12-1-4-1) -- (11-1-2,90.2-4-1+0.6-4-1+0.6);
	\draw (13-1-2,84-4-1) -- (13.5-1-2,84.6-4-1);
	\draw (12-1-2,85.2-4-1) -- (12.5-1-2,85.8-4-1);
	\draw (11.5-1-2,90.2-4-1+0.6-4-1) -- (11-2,90.2-4-1+0.6-4-1+0.6);
	\end{tikzpicture}
	
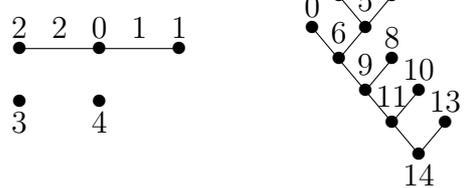
\begin{tikzpicture}[scale=0.7]
	\draw (.3,78) node {};
	\draw (2,77.1) node {$c_i=5$:};
	\draw (5.5-2,98-16-1-5) node {$\bullet$};
	\draw (7-2,98-16-1-5) node {$\bullet$};
	\draw (8.5-2,98-16-1-5) node {$\bullet$};
	\draw (5.5-2,97-16-1-5) node {$\bullet$};
	\draw (7-2,97-16-1-5) node {$\bullet$};
	\draw (5.5-2,98-16-1-5) -- (8.5-2,98-16-1-5);
	\draw (5.5-2,98+0.4-16-1-5) node {2};
	\draw (6.25-2,98.4-16-1-5) node {2};
	\draw (7-2,98+0.4-16-1-5) node {0};
	\draw (7.75-2,98.4-16-1-5) node {1};
	\draw (8.5-2,98+0.4-16-1-5) node {1};
	\draw (5.5-2,97-0.4-16-1-5) node {3};
	\draw (7-2,97-0.4-16-1-5) node {4};
	\draw (13-1-1,89-4-1-4-1-5) node {$\bullet$};
	\draw (13-1-1,88.6-4-1-4-1-5) node {14};
	\draw (12.5-1-1,89.6-4-1-4-1-5) node {$\bullet$};
	\draw (13.5-1-1,96.6-3-4-4-1-4-1-5) node {$\bullet$};
	\draw (12-1-1,90.2-4-1-4-1-5) node {$\bullet$};
	\draw (12-1-1,90.7-4-1-4-1-5) node {9};
	\draw (13-1-1,90.2-4-1-4-1-5) node {$\bullet$};
	\draw (13-1-1,90.6-4-1-4-1-5) node {10};
	\draw (11.5-1-1,90.2-4-1+0.6-4-1-5) node {$\bullet$};
	\draw (11.5-1-1,86.3-4-1-5) node {6};
	\draw (12.5-1-1,85.8-4-1-5) node {$\bullet$};
	\draw (12.5-1-1,86.2-4-1-5) node {8};
	\draw (11-1-1,90.2-4-1+0.6-4-1+0.6-5) node {$\bullet$};
	\draw (11-1-1,90.2-4-1+0.6-4-1+1-5) node {0};
	\draw (12.5-1-1,96.6-3-4-4-1-4-1-5) -- (13-1-1,90.2-4-1-4-1-5);
	\draw (12.5-1-1,97.1-3-4-4-1-4-1-5) node {11};
	\draw (13.5-1-1,97-3-4-4-1-4-1-5) node {13};
	\draw (11-1,90.2-4-1+0.6-4-1+0.6-5) node {$\bullet$};
	\draw (11-1,90.3-4-1+0.6-4-1+1-5) node {5};
	\draw (13-1-1,97-12-1-4-1-5) -- (11-1-1,90.2-4-1+0.6-4-1+0.6-5);
	\draw (13-1-1,84-4-1-5) -- (13.5-1-1,84.6-4-1-5);
	\draw (12-1-1,85.2-4-1-5) -- (12.5-1-1,85.8-4-1-5);
	\draw (11.5-1-1,90.2-4-1+0.6-4-1-5) -- (11-1,90.2-4-1+0.6-4-1+0.6-5);
	\draw (11.5-1,90.2-4-1+0.6-4-1+0.6-5+0.6) node {$\bullet$};
	\draw (10.5-1,90.2-4-1+0.6-4-1+0.6-5+0.6) node {$\bullet$};
	\draw (11.5-1,90.2-4-1+0.6-4-1+0.6-5+1) node {3};
	\draw (10.5-1,90.2-4-1+0.6-4-1+0.6-5+1) node {4};
	\draw (11-1,90.2-4-1+0.6-4-1+0.6-5) -- (11.5-1,90.2-4-1+0.6-4-1+0.6-5+0.6);
	\draw (11-1,90.2-4-1+0.6-4-1+0.6-5) -- (10.5-1,90.2-4-1+0.6-4-1+0.6-5+0.6);
	\draw (13-1,88.6-4-1-4-1-5.5) node {};
	\end{tikzpicture}\\
	There are no more critical edges left, so the construction of $M(X,f)$ is finished.
\end{Ex}
Since there are different notions of merge trees in the literature and since the merge trees in our setting carry a lot of structure, there are multiple possibilities of how to define equivalences of merge trees. In the remainder of this section, we define and discuss some versions of equivalence of merge trees.\par
Since merge trees are defined to be chiral rooted binary trees, the obvious notion for isomorphisms of merge trees is isomorphisms of chiral rooted binary trees. In detail, this means bijections between the sets of nodes and the sets of vertices which map the root to the root and are compatible with the chiral child relation. 
For Mo trees (\Cref{morseorder}) we have more notions of equivalence.
\begin{Def}\label{MoTiso}
Let $(T,\leq)$ and $(T',\leq')$ be Mo trees. An \D{isomorphism} of Mo trees $(T,\leq)\cong (T',\leq')$ is an order-preserving isomorphism of the underlying merge trees. \par
Let $(T,\lambda)$ and $(T',\lambda')$ be Ml trees. An \D{isomorphism} of Ml trees is and isomorphism of the underlying merge trees over $\mathbb{R}$, that is, an isomorphism of merge trees $\varphi\colon T\rightarrow T'$ such that $ \lambda' \circ \varphi=\lambda$.
%a pair $(\varphi,\psi)\colon (T,\lambda) \rightarrow (T',\lambda')$ such that $\varphi \colon T\rightarrow T'$ is an isomorphism of the underlying merge trees, $\psi \colon \mathbb{R}\rightarrow \mathbb{R}$ is a bijection, the restriction $\psi_{\lvert \operatorname{im}(\lambda)}\colon \operatorname{im}(\lambda) \rightarrow \operatorname{im}(\lambda')$ is an order-preserving bijection, and $\psi \circ \lambda = \lambda' \circ \varphi$ holds. 
\end{Def}
\begin{Def}\label{MlTiso}
Let $(T,\lambda)$ and $(T',\lambda')$ be Ml trees. A \D{shuffle equivalence} $(\varphi,\psi)\colon(T,\lambda)\rightarrow (T',\lambda')$ of Ml trees is a pair of an isomorphism of the underlying merge trees $\varphi\colon T\rightarrow T'$ and a bijection $\psi\colon \mathbb{R}\rightarrow\mathbb{R}$ such that
\begin{itemize}
    \item $\psi \circ \lambda = \lambda' \circ \varphi$ holds,
    \item the restriction of $\psi$ to values on leaves is order-preserving, and
    \item the restriction of $\psi$ to values on inner nodes is order-preserving.
\end{itemize}
In the special case that the restriction $\psi_{\lvert \operatorname{im}(\lambda)}\colon \operatorname{im}(\lambda)\rightarrow \operatorname{im}(\lambda')$ is an order preserving bijection, we call $(\varphi,\psi)$ an \D{order equivalence}.\par
A \D{shuffle equivalence} $(T,\leq)\rightarrow (T',\leq')$ between Mo trees is an isomorphism $\varphi$ of the underlying merge trees such that 
\begin{itemize}
    \item the restriction of $\varphi$ to leaf nodes is order-preserving, and
    \item the restriction of $\varphi$ to inner nodes is order-preserving.
\end{itemize}
\end{Def}
\begin{Rem}
The name of shuffle equivalences hints at the fact that two given total orders, one total order on the leaves, and another total order on the inner nodes, might be combined to produce a total order on all nodes in different ways, a bit like shuffling cards. Shuffle equivalence checks if two Morse orders arise from the same underlying orders by different ways of shuffling. But the necessity of ranking ancestors higher than descendants and being compatible with the chirality prevents arbitrary ways of shuffling two given orders on the leaves and inner nodes from producing Morse orders.
\end{Rem}
Shuffle equivalences induce by definition isomorphisms on the underlying merge trees. We now make the relationship between Mo trees and Ml trees precise.
\begin{Prop}\label{moml}
    The Morse labeling induced by a Morse order and the order induced by a Morse labeling define inverse bijections $\operatorname{iMl} \colon MoT/_{\cong} \xleftrightarrow{\cong} MlT/_{\sim}\colon \operatorname{iMo}$ where $\sim$ denotes order equivalence.
\end{Prop}
\begin{proof}
    Let $(T,\leq)$ be a Mo tree. It is immediate that $(T,\lambda_\leq)$ has the property that $\lambda_\leq$ induces $\leq$ as its induced order on $T$. Thus, the composition $\operatorname{iMo}\circ\operatorname{iMl}$ is the identity on $MoT$.\par
    Let $(T,\lambda)$ be an Ml tree. By definition, the labeling $\lambda$ induces a Morse order $\leq_\lambda$ on $T$ which makes $(T,\leq_\lambda)$ a Mo tree. Since the induced Morse labeling $\lambda_{\leq_\lambda}$ by construction induces $\leq_\lambda$ as its induced order, it follows that $(T,\lambda)$ and $(T,\lambda_{\leq_\lambda})$ are order equivalent.
\end{proof}
\begin{Cor}\label{shuffleequiv}
    The induced Morse order $\operatorname{iMo}$ and the induced Morse labeling $\operatorname{iMl}$ induce a bijection between shuffle equivalences of Mo trees and shuffle equivalences of Ml trees.
\end{Cor}
\begin{proof}
    The assignment $\operatorname{iMl}$ maps shuffle equivalences of Mo trees to shuffle equivalences of Ml trees by \Cref{MlTiso}. \par
    Let $\lambda\colon T\rightarrow \mathbb{R} $ and $\lambda' \colon T' \rightarrow \mathbb{R}$ be shuffle-equivalent Morse labelings and let $(\varphi,\psi)$ be the corresponding shuffle equivalence. Then by \Cref{ML} the map $\operatorname{iMo}(\varphi,\psi)=\varphi\colon (T,\leq_\lambda)\rightarrow (T',\leq_{\lambda'})$ has the property that
%    \begin{itemize}
%        \item 
        the restriction of $\varphi$ to leaf nodes is order-preserving, and 
%        \item 
        the restriction of $\varphi$ to inner nodes is order-preserving.
%    \end{itemize}
Hence, $\varphi=\operatorname{iMo}(\varphi,\psi)$ is a shuffle equivalence of Mo trees.
\end{proof}
\begin{Rem}
In particular, the aforementioned proposition and corollary mean that two Mo trees are isomorphic (respectively shuffle equivalent) if and only if the corresponding Ml trees are order equivalent (respectively shuffle equivalent) and vice versa.
\end{Rem}
\subsection{Generic Properties and Equivalences of DMFs}\label{slc}
In this subsection, we will take a closer look at generic properties that dMfs can be assumed to have. Furthermore we consider some notions of equivalences between dMfs.\par 
A first example of a generic property is the notion of index-ordered dMfs as defined above. It is inspired by the eponymous notion from the smooth case. We use index-ordered dMfs in order to distinguish critical simplices by their dimension because merge trees have the same property: Critical 0-simplices appear as leaves whereas critical 1-simplices appear as inner nodes of the induced merge tree (see \Cref{Mer}). Thus, index-ordered dMfs seem to be especially suitable for working with merge trees. \par
Nonetheless, index-ordered dMfs are not compatible with the structure of rooted subtrees. In detail, consider the following:
\begin{Rem}\label{muex}
Let $f\colon X \rightarrow \mathbb{R}$ be an index-ordered dMf on a tree such that the following holds: The tree $X$ has a critical 1-simplex $\sigma$ such that the corresponding inner node $p$ in $M(X,f)$ has two inner nodes $c$ and $c'$ as children. Then the image of $f$ on at least one of the connected components corresponding to $c$ or $c'$ is not an interval in $f(\operatorname{Cr}(f))$.
\end{Rem}
\begin{Ex}
	The following is a small example:\par
	
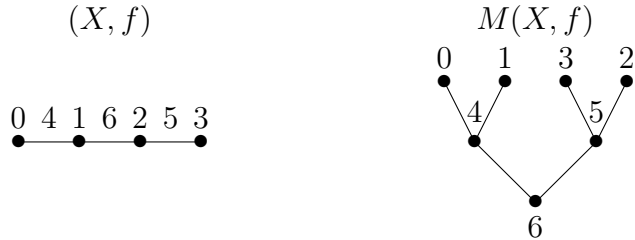
\begin{tikzpicture}[scale=0.8]
	\draw (-2,3.3) node {};
	\draw (4,3) node {$(X,f)$};
	\draw (2.5,1) node {$\bullet$};
	\draw (3.5,1) node {$\bullet$};
	\draw (4.5,1) node {$\bullet$};
	\draw (5.5,1) node {$\bullet$};
	\draw (2.5,1.4) node {$0$};
	\draw (3.5,1.4) node {$1$};
	\draw (4.5,1.4) node {$2$};
	\draw (5.5,1.4) node {$3$};
	\draw (2.5,1) -- (5.5,1);
	\draw (3,1.4) node {$4$};
	\draw (4,1.4) node {$6$};
	\draw (5,1.4) node {$5$};
	\draw (11,3) node {$M(X,f)$};
	\draw (11,0) node {$\bullet$};
	\draw (11,-0.4) node {$6$};
	\draw (11,0) -- (10,1);
	\draw (11,0) -- (12,1);
	\draw (10,1) node {$\bullet$};
	\draw (10,1.5) node {$4$};
	\draw (10,1) -- (9.5,2);
	\draw (10,1) -- (10.5,2);
	\draw (12,1) node {$\bullet$};
	\draw (12,1.5) node {$5$};
	\draw (12,1) -- (11.5,2);
	\draw (12,1) -- (12.5,2);
	\draw (9.5,2) node {$\bullet$};
	\draw (10.5,2) node {$\bullet$};
	\draw (11.5,2) node {$\bullet$};
	\draw (12.5,2) node {$\bullet$};
	\draw (9.5,2.4) node {$0$};
	\draw (10.5,2.4) node {$1$};
	\draw (11.5,2.4) node {$3$};
	\draw (12.5,2.4) node {$2$};
	\end{tikzpicture}\par
	Neither the subtree with root labeled $4$, nor the subtree with root labeled $5$ is labeled with an interval in $f(\operatorname{Cr}(f))\subseteq\mathbb{R}$.
\end{Ex} 
Still, we can assume compatibility with the structure of rooted subtrees as a property:
\begin{Def}\label{defsubcon}
	Let $f\colon X \rightarrow \mathbb{R}$ be a dMf. The function $f$ is called \D{sublevel-connected} if for all critical 1-simplices $v$ the set $f(X_{f(v)}[v])$ is an interval in $f(\operatorname{Cr}(f))$.
\end{Def}
\begin{Rem}
	Since both \lq{}index-ordered\rq{} and \lq{}sublevel-connected\rq{} are properties that only rely on the values of $f$ on critical simplices, they can easily be arranged without changing the partial matching if the tree $X$ is finite. It is also possible to do this without changing the induced merge tree. One way to obtain a sublevel-connected dMf would be to choose a collapsing order for the induced matching such that the respective connected components of sublevel sets correspond to intervals in said collapsing order. But as seen in \Cref{muex}, the two properties are in most cases mutually exclusive.
\end{Rem}
\begin{Ex}
	Here is a possibility how to modify the dMf from the previous example in order to make it sublevel-connected without changing the induced merge tree:\par
	
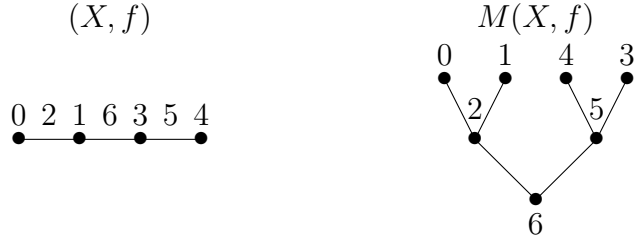
\begin{tikzpicture}[scale=0.8]
	\draw (-2,1) node {};
	\draw (4,3) node {$(X,f)$};
	\draw (2.5,1) node {$\bullet$};
	\draw (3.5,1) node {$\bullet$};
	\draw (4.5,1) node {$\bullet$};
	\draw (5.5,1) node {$\bullet$};
	\draw (2.5,1.4) node {$0$};
	\draw (3.5,1.4) node {$1$};
	\draw (4.5,1.4) node {$3$};
	\draw (5.5,1.4) node {$4$};
	\draw (2.5,1) -- (5.5,1);
	\draw (3,1.4) node {$2$};
	\draw (4,1.4) node {$6$};
	\draw (5,1.4) node {$5$};
	\draw (11,3) node {$M(X,f)$};
	\draw (11,0) node {$\bullet$};
	\draw (11,-0.4) node {$6$};
	\draw (11,0) -- (10,1);
	\draw (11,0) -- (12,1);
	\draw (10,1) node {$\bullet$};
	\draw (10,1.5) node {$2$};
	\draw (10,1) -- (9.5,2);
	\draw (10,1) -- (10.5,2);
	\draw (12,1) node {$\bullet$};
	\draw (12,1.5) node {$5$};
	\draw (12,1) -- (11.5,2);
	\draw (12,1) -- (12.5,2);
	\draw (9.5,2) node {$\bullet$};
	\draw (10.5,2) node {$\bullet$};
	\draw (11.5,2) node {$\bullet$};
	\draw (12.5,2) node {$\bullet$};
	\draw (9.5,2.4) node {$0$};
	\draw (10.5,2.4) node {$1$};
	\draw (11.5,2.4) node {$4$};
	\draw (12.5,2.4) node {$3$};
	\end{tikzpicture}
\end{Ex}
Since dMfs can be modified to fulfill either of the two properties, one can always choose the one which is more convenient for the task at hand. Thus, we give two different constructions in this work, one for each property. \par
We recall that by \Cref{critorder} each dMf $f \colon X \rightarrow \mathbb{R}$ induces an order on the 0-simplices of $X$ and on the 1-simplices of $X$, respectively.
This yields the following definition, which induces a merge-tree-invariant notion of equivalence between dMfs.
\begin{Def}
     Let $f\colon X\rightarrow \mathbb{R}$ and $g\colon X\rightarrow \mathbb{R}$ be dMfs on a tree $X$. We call $f$ and $g$ \D{shuffle-equivalent} if they have the same critical simplices and if they induce the same order on the critical 0-simplices as well as the same order on the critical 1-simplices. \par
     Let $(X,f)$ and $(X',f')$ be two dMfs on trees. A \D{shuffle equivalence} $(\varphi,\psi)\colon f \rightarrow f'$ between $f$ and $f'$ consists of a simplicial map $\varphi\colon X \rightarrow X'$ and a bijection $\psi\colon \mathbb{R} \rightarrow \mathbb{R}$ such that 
     \begin{itemize}
         \item $\psi \circ f=f' \circ \varphi$, 
         \item $\varphi_{\lvert \operatorname{Cr}(f)}\colon \operatorname{Cr}(f) \rightarrow \operatorname{Cr}(f')$ is a bijection,
         \item the restriction of $\psi$ to values on critical 0-simplices is order preserving, and
         \item the restriction of $\psi$ to values on critical 1-simplices is order preserving.
     \end{itemize}
 In the special case that the restriction $\psi_{\lvert \operatorname{Cr}(f)}\colon \operatorname{Cr}(f)\rightarrow \operatorname{Cr}(f')$ is an order preserving bijection, we call $(\varphi,\psi)$ an \D{order equivalence}.\par
\end{Def}
\begin{Rem}
    It is immediate that shuffle equivalence of dMfs is an equivalence relation. The name is inspired analogously as for the eponymous notion for Mo trees and Ml trees.
    Shuffle equivalence checks if two dMfs arise from the same underlying orders by different ways of shuffling. Nonetheless not all ways of shuffling given orders on the critical 0-simplices and critical 1-simplices produce a dMf because dMfs have to be weakly increasing.\par
    Furthermore, shuffle equivalence, and in particular order equivalence, also considers dMfs to be equivalent if they only differ by scaling because a different scaling does not change the induced orders on simplices.
\end{Rem}
We split the definition of shuffle equivalences in two steps to simplify the proofs of the following propositions. It is immediate that for two dMfs on the same tree $X$ there is a shuffle equivalence between them if and only if they are shuffle equivalent.
\begin{Prop}\label{equiv}
Let $X$ be a tree and let $f\colon X \rightarrow \mathbb{R}$ and $g\colon X \rightarrow \mathbb{R}$ be two shuffle-equivalent dMfs.
Then $M(X,f)$ and $M(X,g)$ are isomorphic as merge trees.
\end{Prop}
The following two lemmas will be helpful for the proof of the proposition:\par
\begin{Lem}\label{critmax}
    Let $X$ be a tree and let $f\colon X \rightarrow \mathbb{R}$ be a dMf such that $X$ has at least one critical 1-simplex. Then the function $f_{\lvert \operatorname{Cr}(f)}$ attains its maximum on a critical 1-simplex. 
\end{Lem}
\begin{proof}[Sketch of Proof]Backtracking gradient paths, see \cite[Def 8.4]{Mor}, as long as possible leads to a local maximum which turns out to be a critical 1-simplex. 
%    We prove the lemma by contradiction. Assume that $f_{\lvert \operatorname{Cr}(f)}$ attains its maximum on a critical 0-simplex $\sigma$. Let $\tau$ be a critical 1-simplex of $X$ such that there is no other critical 1-simplex on the unique path between $\sigma$ and $\tau$. Since $X$ is a tree, such a unique path always exists. Said unique path consists of a chain of simplices $\sigma \subset a_0 \supset a_1 \subset \dots \supset a_n \subset \tau$ such that all simplices $a_i$ are non-critical. Then the sequence induces the sequence $f(\sigma)<f(a_0) =f(a_1) < \dots < f(a_n) < f(\tau)$. This is a contradiction to the assumption that $f_{\lvert \operatorname{Cr}(f)}$ attains its maximum on $\sigma$.
\end{proof}
%\begin{Rem}
%    The chain of simplices used in the previous proof is an example of a gradient path, see \cite[Def 8.4]{Mor}. 
%\end{Rem}
\begin{Lem}\label{subcon}
    Let $X$ be a tree and let $f \colon X \rightarrow \mathbb{R}$ and $g \colon X \rightarrow \mathbb{R}$ be two shuffle-equivalent dMfs.
    We denote the critical 1-simplices of $f$ and $g$ by $c_0<c_1< \dots <c_n$ where $<$ denotes the ordering induced by $f$ or $g$, respectively. Then the connected components $X^f_{f(c_i)}[c_j]$ of sublevel complexes contain the same critical simplices as $X^g_{g(c_i)}[c_j]$ for all $j<i$.
\end{Lem}
\begin{proof}[Sketch of Proof]
	First we observe that restrictions of dMfs on trees to connected components of sublevel complexes are again dMfs on trees. With help of \Cref{critmax}, it can be proved inductively that in the construction of $M(X,f)$ and $M(X,g)$ the same critical 1-simplices are considered in the same order. The statement then follows inductively. 
\end{proof}
%\begin{proof} 
%    Connected components of sublevel complexes are trees themselves. Furthermore, the restrictions of dMfs to subcomplexes are again dMfs because being a dMf only depends on properties that are inherited by arbitrary subcomplexes of $X$. It follows by \Cref{critmax} that for connected components $X^f_{f(c_i)}[c_j]/X^g_{g(c_i)}[c_j]$ which contain a critical 1-simplex $c_j$, the maps $f_{\lvert X^f_{f(c_i)}[c_j]}$  and $g_{\lvert X^g_{g(c_i)}[c_j]}$ attain their maxima on critical 1-simplices. It follows inductively that $f_{\lvert X^f_{f(c_i)}[c_j]}$  and $g_{\lvert X^g_{g(c_i)}[c_j]}$ attain their maxima in each step on the same critical 1-simplex, namely $c_j$, because both functions induce the same order on the critical 1-simplices. Thus, in each step in which a lower level is considered, the connected components of sublevel complexes are divided into smaller pieces at the same critical 1-simplex. It might happen that the connected components of sublevel complexes of lower levels differ by matched 1-simplices that have higher values than the next maximal critical value. However, by \cite[Lem 2.6]{Mor} the homotopy type and, therefore, the connected components of sublevel complexes only change at critical levels. Because of that it follows inductively that critical simplices are located inside the same new connected components of sublevel sets of lower levels for both functions $f$ and $g$. \par 
%    
%\end{proof}
\begin{proof}[Proof of \Cref{equiv}]
We consider the construction of the induced merge tree (see \Cref{Mer}) and prove inductively the slightly stronger result that both functions yield isomorphic merge trees at every step of the construction.
%that for both functions every step yields the same results regarding the isomorphism type of the induced merge tree. 
This implies that $f$ and $g$ induce isomorphic merge trees.\par
Since $f$ and $g$ impose the same order on the set of critical 1-simplices, the construction algorithm considers the same critical 1-simplices during the same steps for both functions. This already proves the base case. In particular, this means that the created root node corresponds to the same critical 1-simplex for both functions. Although the label of the root node might be different for the two dMfs, it does not affect the isomorphism type of the induced merge tree because the labeling is not part of the data of merge trees.  \par
For the inductive step, we observe that in every step of the construction we consider a connected component of sublevel complexes $X^f_{f(c_i)}[c_j]/X^g_{g(c_i)}[c_j]$ that contains at least one critical 1-simplex, namely the one with the highest remaining critical value $c_j$. Thus, by \Cref{subcon} in each step of the construction, the same critical simplices occur. \par
Assume we are at the step that considers the critical 1-simplex $c_i$. For the two new nodes which are created in the inductive step, two pieces of information are important for the isomorphism type of the induced merge tree, namely the chirality of the new nodes and which critical simplices the new nodes correspond to. The chirality of the new nodes affects the isomorphism type of the induced merge tree directly. The critical simplex corresponding to a child node $c$ decides which connected component of the respective sublevel complex is used to build the subtree with root $c$ and at which point said connected component will be subdivided next.\par
Both pieces of information are defined by the two connected components that belong to the boundary 0-simplices of $c_i$. The two child nodes correspond to the critical simplices with the highest critical values. \par
There are three cases:
\begin{enumerate}
    \item Both connected components contain at least one critical 1-simplex $c_j$. 
    \item One connected component contains at least one critical 1-simplex $c_j$ whereas the other one only contains one critical 0-simplex $c$.
    \item Each of the two connected components contains only one critical 0-simplex $c$.
\end{enumerate}
It follows by \Cref{critmax} that in case (1) the corresponding 1-simplices with the highest critical values are critical 1-simplices $c_j$. In case (2) the same is true for the connected component that contains at least one critical 1-simplex. For the connected components in case (2) and (3) that only contain one critical 0-simplex, respectively, it is true that the critical 0-simplex is the only critical simplex in its corresponding connected component. Thus, the new nodes correspond to the same critical simplices for $f$ and for $g$ because both functions induce the same order on 1-simplices and because connected components only correspond to critical 0-simplices if they are the only critical simplex left in the corresponding connected component. Furthermore, connected components $X^f_{f(c_i)}[c]/X^g_{g(c_i)}[c]$ that only contain one critical 0-simplex do not have any influence on the induced merge tree $M(X,f)/M(X,g)$ because their corresponding nodes have already been created during a step that considered a critical 1-simplex with a higher critical value and they are not considered in later steps of the construction.\par
The chirality of the new nodes depends on the minimal values on critical simplices of the two respective connected components $X^f_{f(c_i)}[c_j]/X^g_{g(c_i)}[c_j]$ or $X^f_{f(c_i)}[c]/X^g_{g(c_i)}[c]$. By \Cref{min}, these minima belong to critical 0-simplices. By assumption, $f$ and $g$ induce the same order on the critical 0-simplices, so the same 0-simplex is minimal with respect to both functions. Thus, $f$ and $g$ assign the same chirality to the new nodes.

\end{proof}
\begin{Rem}
The functions $f$ and $g$ might induce different order relations between 0-simplices and 1-simplices. Therefore, in sublevel complexes that appear during the same step of the construction there might be different connected components that contain only one critical 0-simplex each with respect to the two dMfs. However, those connected components that contain only one critical 0-simplex and no critical 1-simplices do not affect the isomorphism type of the induced merge tree. This is because those connected components correspond to leaves of the merge tree which have already been created during the step that considered the critical 1-simplex between said connected components and other connected components. Furthermore, those connected components do not appear in later steps of the construction of the induced merge tree because they do not contain any critical 1-simplices.     
\end{Rem}
\begin{Prop}\label{shuffleorder}
    Shuffle equivalences of dMfs induce shuffle equivalences of the induced Ml trees. Moreover, order equivalences of dMfs induce order equivalences of the induced Ml trees.
\end{Prop}
\begin{proof}[Sketch of Proof]
	Since $\varphi$ is bijective on critical simplices and simplicial, it follows that $\varphi$ induces a bijection between connected components of sublevel complexes. With this, it follows analogously to the proof of \Cref{equiv} that $M(X,f)\cong M(X',f')$ holds. The proof that the induced Morse labelings are shuffle equivalent is straightforward and only uses that the restrictions of $\psi$ to 0-simplices and to 1-simplices are order preserving, and the compatibility between $\varphi,\psi,f$ and $f'$.  
\end{proof}
\begin{Rem}
    The criterion from \Cref{equiv} for merge tree equivalence is sufficient but not necessary, as the following example shows:\par
    
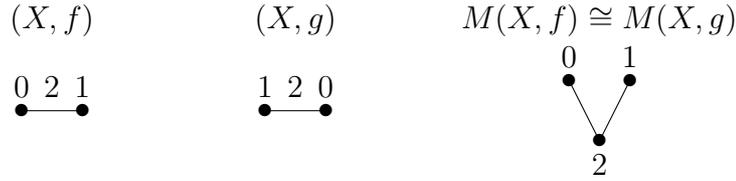
\begin{tikzpicture}[scale=0.8]
        \draw (-1,2.3) node {};
        \draw (3,2) node {$(X,f)$};
        \draw (2.5,0.5) node {$\bullet$};
        \draw (2.5,0.9) node {$0$};
        \draw (3.5,0.5) node {$\bullet$};
        \draw (3.5,0.9) node {$1$};
        \draw (2.5,0.5) -- (3.5,0.5);
        \draw (3,0.9) node {$2$};
        \draw (7,2) node {$(X,g)$};
        \draw (6.5,0.5) node {$\bullet$};
        \draw (6.5,0.9) node {$1$};
        \draw (7.5,0.5) node {$\bullet$};
        \draw (7.5,0.9) node {$0$};
        \draw (6.5,0.5) -- (7.5,0.5);
        \draw (7,0.9) node {$2$};
        \draw (12,2) node {$M(X,f)\cong M(X,g)$};
        \draw (12,0) node {$\bullet$};
        \draw (12,-0.4) node {$2$};
        \draw (11.5,1) node {$\bullet$};
        \draw (11.5,1.4) node {$0$};
        \draw (12,0) -- (11.5,1);
        \draw (12.5,1) node {$\bullet$};
        \draw (12.5,1.4) node {$1$};
        \draw (12,0) -- (12.5,1);
    \end{tikzpicture}\par
    The two dMfs $f$ and $g$ induce inverse orders on the two 0-simplices. Nonetheless, $f$ and $g$ induce the same unlabeled merge tree.
\end{Rem}
This remark leads us to yet another kind of equivalence relation between dMfs that arises from symmetries of sublevel complexes. In order to make this notion of symmetry precise, we need some preparations.
\begin{Def}\label{symm}
     Let $f\colon X \rightarrow \mathbb{R}$ be a dMf on a tree. For each non-empty connected component $X^f_c[v]$ of a sublevel complex $X^f_c$ we denote by $\operatorname{Aut}(X^f_c[v])$ the group of simplicial automorphisms of $X^f_c[v]$. For each $a \in \operatorname{Aut}(X^f_c[v])$ there is an extension to a self-bijection $X\rightarrow X$ by the identity. The group $\widetilde{\operatorname{Aut}}(X^f_c[v])$ is defined to be the group of said extensions of elements of $\operatorname{Aut}(X^f_c[v])$ by the identity. We consider $\widetilde{\operatorname{Aut}}(X^f_c[v])$ as a subgroup of the group of all self-bijections of $X$. The total order on $\operatorname{Cr}(f)$ induced by $f$ induces chains  $\widetilde{\operatorname{Aut}}(X^f_{c_0}[v])\subset \widetilde{\operatorname{Aut}}(X^f_{c_1}[v]) \subset \dots$ of inclusions of subgroups. Moreover, we have inclusions $\widetilde{\operatorname{Aut}}(X^f_{c_i}[v])\subset \widetilde{\operatorname{Aut}}(X^f_{c_j}[v])=\widetilde{\operatorname{Aut}}(X^f_{c_j}[w]) \supset \widetilde{\operatorname{Aut}}(X^f_{c_i}[w])$ if $v$ and $w$ are in different connected components of some sublevel complex $X^f_{c_i}$ that merge together in some other sublevel complex $X^f_{c_j}$ for $j>i$.  We define the \D{sublevel automorphism group} of $(X,f)$, denoted by $\operatorname{Aut}_{sl}(X,f)$, to be the subgroup generated by $ \bigcup\limits_{c \in \operatorname{Cr}(f),v \in X} \widetilde{\operatorname{Aut}}(X^f_c[v])$. We call the elements of $\operatorname{Aut}_{sl}(X,f)$ \D{sublevel automorphisms}.
\end{Def}
\begin{Rem}
	Even though sublevel automorphisms are built out of simplicial automorphisms of connected components $X^f_c[v]$ of sublevel complexes, they are in general not simplicial maps $X \rightarrow X$. To be precise, if a simplicial automorphism of $X^f_c[v]$ was used to construct a sublevel automorphism $a\in \operatorname{Aut}_{sl}(X,f)$, then $a$ will fail to be simplicial at the boundary of $X^f_c[v]\subset X$.
\end{Rem}
\begin{Prop}\label{autgrpac}
    Let $f\colon X\rightarrow \mathbb{R}$ be a dMf on a tree and let $a \in \operatorname{Aut}_{sl}(X,f)$ be a sublevel automorphism.
    Then $f*a$ defined by $f*a(\sigma)\coloneqq f(a(\sigma))$ is a dMf on $X$. Moreover, this defines a right group action of $\operatorname{Aut}_{sl}(X,f)$ on the set of dMfs on $X$.
\end{Prop}
\begin{proof}[Sketch of Proof]
	The proof that $f*a$ is a dMf is straightforward and the compatibility of the group action follows directly by associativity of the composition of maps.
\end{proof}
%\begin{proof}
%    Since $a$ is in particular a self-bijection of $X$, the map $f*a$ is distinct-valued and at most 2-1 because $f$ is. Because $a$ is an automorphism of a sublevel complex $X^f_c$, it preserves the face relations between cells of $X^f_c$. Thus, $f*a$ is weakly increasing on $X^f_c$ because $f$ is. The function $f*a$ is also weakly increasing at the boundary of $X_c^f$ because $X_c^f$ is a sublevel complex. \par
%    The compatibility of the group action follows directly by the associativity of the composition of maps.
%\end{proof}
\begin{Rem}
Since automorphisms of simplicial complexes preserve the dimension of simplices, the action of $\operatorname{Aut}_{sl}(X,f)$ on the set of dMfs on a tree $X$ preserves the properties of being index-ordered or sublevel-connected. 
\end{Rem}
\begin{Def}
     Let $f\colon X\rightarrow \mathbb{R}$ and $g\colon X\rightarrow \mathbb{R}$ be dMfs on a tree $X$. We call $f$ and $g$ \D{sublevel-equivalent} if $\operatorname{Cr}(f)=\operatorname{Cr}(g)$ and $X^f_c\cong X^g_c$ for all $c\in \operatorname{Cr}(f)=\operatorname{Cr}(g)$. If additionally $g=f*a$ holds for a sublevel automorphism $a \in \operatorname{Aut}_{sl}(X,f)=\operatorname{Aut}_{sl}(X,g)$, then we call $f$ and $g$ \D{symmetry-equivalent}. We call the map $a$ a \D{symmetry equivalence} from $f$ to $g$. \par
     We call two dMfs $f\colon X\rightarrow \mathbb{R}$ and $g\colon Y\rightarrow \mathbb{R}$ \D{symmetry-equivalent} if there is a simplicial isomorphism $\varphi\colon X\rightarrow Y$ such that $f$ and $g\circ \varphi$ are symmetry-equivalent.
\end{Def}
\begin{Ex}\label{symmequivalentdmf}
    We give a list of some symmetry-equivalent dMfs on a path with four vertices. Here we denote the sublevel equivalence induced by the reflection of the connected component of the sublevel complex of level $k$ by $a_k[k]$ :\\
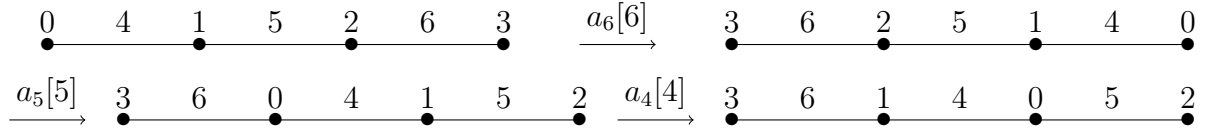
    \begin{tikzpicture}
    \draw (0,9.5) node {};
    \draw (0,9) node {$\bullet$};
    \draw (2,9) node {$\bullet$};
    \draw (4,9) node {$\bullet$};
    \draw (6,9) node {$\bullet$};
    \draw (0,9) -- (6,9);
    \draw (0,9.3) node {$0$};
    \draw (1,9.3) node {$4$};
    \draw (2,9.3) node {$1$};
    \draw (3,9.3) node {$5$};
    \draw (4,9.3) node {$2$};
    \draw (5,9.3) node {$6$};
    \draw (6,9.3) node {$3$};
    \draw[->] (7,9) -- (8,9) node[midway, above] {$a_6[6]$};
    \draw (9,9) node {$\bullet$};
    \draw (11,9) node {$\bullet$};
    \draw (13,9) node {$\bullet$};
    \draw (15,9) node {$\bullet$};
    \draw (9,9) -- (15,9);
    \draw (9,9.3) node {$3$};
    \draw (10,9.3) node {$6$};
    \draw (11,9.3) node {$2$};
    \draw (12,9.3) node {$5$};
    \draw (13,9.3) node {$1$};
    \draw (14,9.3) node {$4$};
    \draw (15,9.3) node {$0$};
    \draw[->] (-.5,8) -- (.5,8) node[midway, above] {$a_5[5]$};
    \draw (0+1,8) node {$\bullet$};
    \draw (2+1,8) node {$\bullet$};
    \draw (4+1,8) node {$\bullet$};
    \draw (6+1,8) node {$\bullet$};
    \draw (0+1,8) -- (6+1,8);
    \draw (0+1,8.3) node {$3$};
    \draw (1+1,8.3) node {$6$};
    \draw (2+1,8.3) node {$0$};
    \draw (3+1,8.3) node {$4$};
    \draw (4+1,8.3) node {$1$};
    \draw (5+1,8.3) node {$5$};
    \draw (6+1,8.3) node {$2$};
    \draw[->] (7.5,8) -- (8.5,8) node[midway, above] {$a_4[4]$};
    \draw (9,8) node {$\bullet$};
    \draw (11,8) node {$\bullet$};
    \draw (13,8) node {$\bullet$};
    \draw (15,8) node {$\bullet$};
    \draw (9,8) -- (15,8);
    \draw (9,8.3) node {$3$};
    \draw (10,8.3) node {$6$};
    \draw (11,8.3) node {$1$};
    \draw (12,8.3) node {$4$};
    \draw (13,8.3) node {$0$};
    \draw (14,8.3) node {$5$};
    \draw (15,8.3) node {$2$};
    \end{tikzpicture}
\end{Ex}
\begin{Prop}\label{symequiv}
    Let $f\colon X\rightarrow \mathbb{R}$ and $g\colon X\rightarrow \mathbb{R}$ be symmetry-equivalent dMfs on a tree $X$. Then $M(X,f)$ and $M(X,g)$ are isomorphic as Ml trees.
\end{Prop}
\begin{proof}[Sketch of Proof]
	The proposition can be proved by induction over the level $c$ of the sublevel automorphisms that the given symmetry equivalence consists of. We check that in each step the single sublevel automorphism $a$ of the connected component $X_c[\sigma]$ with the 1-simplex $\sigma$ labeled $c$ in $X_c$ only affects steps of the construction of the induced Ml tree that consider simplices of $X_c[\sigma]$. Moreover, it is straightforward to prove that in these steps, the created nodes and their induced Morse labels are the same as without the application of $a$.
\end{proof}
%\begin{proof}
%It suffices to consider the case of a single sublevel automorphism $a\in \widetilde{\operatorname{Aut}}(X_c^f[v])$. Since symmetry equivalences consist of compositions of sublevel automorphisms, the statement then follows by repeated application of the aforementioned case.\par 
%We consider how the application of a sublevel automorphism $a$ affects the steps of the construction of the induced Ml tree, \Cref{Mer}.
%The sublevel automorphism $a$ does not affect steps that do not remove any critical edge of $X_c^f[v]$ because $a$, by definition, acts identically outside of $X_c^f[v]$. Considering a critical edge of $X_c^f[v]$ in the constructing algorithm of the induced Ml tree corresponds to passing from a sublevel complex $X^f_c$ in the sublevel filtration to a sublevel complex of a slightly lower level $X^f_{c-\varepsilon}$. 
%%When a critical edge of $X_c^f[v]$ is being considered, 
%We observe that $a$ is precisely a simplicial automorphism of $X_c^f[v]$. Since simplicial automorphisms preserve the filtration of the corresponding sublevel complex $X_c^f$ into further sublevel complexes $X_{c-\varepsilon}^f$, the parts that $X_c^f[v]$ is divided into by passing to a slightly lower level $c-\varepsilon$ induce the same child nodes with the same chirality and label in the induced Ml tree.\par
%In conclusion, the construction yields the same results as without the application of $a$ during all steps. 
%\end{proof}
\begin{Rem}
Sublevel automorphisms of dMfs on paths only consist of reflections of the corresponding connected component of a subcomplex. 
When such a connected component of a sublevel complex of a critical level $c$ is considered during the construction of the induced Ml tree, the reflection of the connected component only causes the two new parts that are obtained by considering a slightly lower level $c-\varepsilon$ to appear as their mirror images. 
In particular, the given dMf attains the same values on the two new parts as before. Hence, the two parts that appear are simplicially isomorphic to the ones that appear without application of the reflection.  
\end{Rem}

\begin{Def}\label{cmequiv}
	Let $(X,f)$ and $(X',f')$ be dMfs on trees. A \D{component-merge equivalence (cm equivalence) of level $a$} is a bijection $\varphi \colon X \rightarrow X'$ such that one of the following, not necessarily exclusive, cases holds:
	\begin{enumerate}
		\item $\varphi$ is a symmetry equivalence.
		\item $\varphi$ fulfills the following:
	\begin{itemize}
		\item $f' \circ \varphi=f$,
		\item $\varphi$ induces a bijection between the sets of connected components of sublevel complexes such that each restriction $\varphi_{\lvert X_{a-\varepsilon}[v]}\colon X_{a-\varepsilon}[v] \rightarrow X'_{a-\varepsilon}[\varphi(v)]$ is a cm equivalence of some level $b<a$, and
		\item the edge $\sigma \in X$ with $f(\sigma)=a$ merges the two connected components $X_{a-\varepsilon}[v_1]$ and $X_{a-\varepsilon}[v_2]$ in $X_{a}[v_1]=X_{a}[v_2]$ if and only if the edge $\varphi(\sigma)$ merges the two connected components $X'_{a-\varepsilon}[\varphi(v_1)]$ and $X'_{a-\varepsilon}[\varphi(v_2)]$ in $X'_{a}[\varphi(v_1)]=X'_{a}[\varphi(v_2)]$.
	\end{itemize}
If $\varphi$ fulfills property (2) but not property (1), we call $\varphi$ \D{non-trivial}.
\end{enumerate}
\end{Def}
\begin{Ex}\label{cmequivex}
	We give an example of two cm-equivalent dMfs on trees. The non-trivial cm equivalence from the left-hand-side to the right-hand-side consists of a symmetry equivalence of level $5$ and the attachment of the edge labeled $6$ between the vertices labeled 1 and 3 rather than 2 and 3. That is, it is a cm-equivalence of level 6.\\
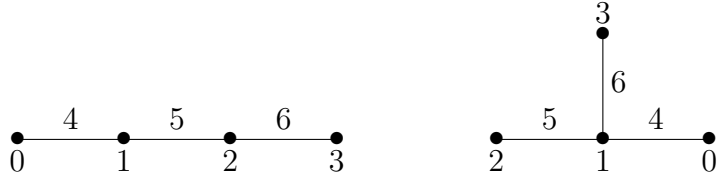
	\begin{tikzpicture}[scale=0.7]
	\draw (-3,2.5) node {};
	\draw (1,0) node {$\bullet$};
	\draw (3,0) node {$\bullet$};
	\draw (5,0) node {$\bullet$};
	\draw (7,0) node {$\bullet$};
	\draw (1,-0.4) node {$0$};
	\draw (2,0.4) node {$4$};
	\draw (3,-0.4) node {$1$};
	\draw (4,0.4) node {$5$};
	\draw (5,-0.4) node {$2$};
	\draw (6,0.4) node {$6$};
	\draw (7,-0.4) node {$3$};
	\draw (1,0) -- (7,0);
	\draw (10,0) node {$\bullet$};
	\draw (12,0) node {$\bullet$};
	\draw (14,0) node {$\bullet$};
	\draw (12,2) node {$\bullet$};
	\draw (10,-0.4) node {$2$};
	\draw (12,-0.4) node {$1$};
	\draw (14,-0.4) node {$0$};
	\draw (12,2.4) node {$3$};
	\draw (11,0.4) node {$5$};
	\draw (10,0) -- (14,0);
	\draw (13,0.4) node {$4$};
	\draw (12,0) -- (12,2);
	\draw (12.3,1.1) node {$6$};
	\end{tikzpicture}
\end{Ex}
\begin{Prop}\label{cmequiviso}
	Cm equivalent dMfs on trees induce isomorphic Ml trees.
\end{Prop}
\begin{proof}
	Let $\varphi \colon (X,f) \rightarrow (X',f')$ be a cm equivalence. By property (ii) of \Cref{dmf}, at most one non-trivial cm equivalence of level $a$ can occur for any level $a$ because there is at most one edge labeled $a$ in $(X,f),(X',f')$, respectively. Thus, we can decompose any cm equivalence into a sequence $(\varphi_a)_a$ of non-trivial cm equivalences of decreasing levels such that each $\varphi_a$ only changes the attachment of the single edge $\sigma$ with $f(\sigma)=a$ and acts as a symmetry equivalence on the rest of path and dMf. It suffices to consider a single level $a$ because the statement then follows by induction from highest to lowest over all levels $a$. \par
	For such a non-trivial cm equivalence $\varphi_a$ we consider the step of the construction of the induced Ml trees that considers the edge $\sigma$ with $f(\sigma)=a$ and the edge $\varphi(\sigma)$. We inductively assume that $\varphi$ induces an isomorphism of induced Ml trees everywhere outside the subtrees corresponding to the two connected components of $X^f_{a-\varepsilon}$ that are merged by the edge $\sigma$ with $f(\sigma)=a$. That is, on the rest of $M(X,f)$ the map $M(\varphi)$ is a bijection compatible with the chiral child relation onto $M(X',f')$ except for the subtrees of $M(X',f')$ which correspond to the connected components of $X'^{f'}_{a-\varepsilon}$ which are merged by the edge $\varphi(\sigma)$. Since the map $\varphi$ is compatible with the dMfs and because it restricts to a cm equivalence $X^f_{a-\varepsilon}\rightarrow X'^{f'}_{a-\varepsilon}$, the dMf $f$ attains the same minima and maxima on the two relevant connected components of $X^{f}_{a-\varepsilon}$ as $f'$ does on their counterparts of $X'^{f'}_{a-\varepsilon}$ via $\varphi$. Since \Cref{Mer} only considers which two connected components are merged by the considered edge, it makes no difference for the isomorphism type of the induced Ml trees that $\sigma$ in general merges the two connected components of $X^f_{a-\varepsilon}$ at vertices that do not correspond via $\varphi$ to the ones adjacent to $\varphi(\sigma)$ in $X'^{f'}_{a-\varepsilon}$. Thus, the construction of the induced Ml tree produces nodes with the same chirality and label for both induced Ml trees in the steps that consider $\sigma,\varphi(\sigma)$, respectively. By assumption, the restriction $\varphi_{X^{f}_{a-\varepsilon}}\colon X^f_{a-\varepsilon}\rightarrow X'^{f'}_{a-\varepsilon}$ is a symmetry equivalence, so the isomorphism of Ml trees extends to the subtrees that correspond to the respective connected components. 
\end{proof}
\begin{Prop}\label{cmequivpath}
	Let $(X,f)$ be a dMf on a tree. There is a dMf on a path $(P,f')$ such that $(X,f)$ is cm-equivalent to $(P,f')$.
\end{Prop}
\begin{proof}[Sketch of Proof] A suitable cm equivalence can be constructed inductively by re-attaching 1-simplices of level $a$ that would become the third 1-simplex incident to some 0-simplex in $X_a$. 
%	Let $\sigma$ be the minimal edge of $X$ with respect to $f$ such that $\sigma$ is adjacent to a vertex $v$ which is adjacent to more than 2 edges in $X_{f(\sigma)}$. Then the connected component $X_{f(\sigma)-\varepsilon}[v]$ is a path and the vertex $v$ is adjacent to exactly two edges in $X_{f(\sigma)-\varepsilon}[v]$. We construct a cm equivalence of level $f(\sigma)$ by reattaching $\sigma$ two one of the two vertices of $X_{f(\sigma)-\varepsilon}[v]$ that are only incident to one edge, respectively. If the other boundary vertex $w$ of $\sigma$ is incident to more than two edges in $X_{f(\sigma)-\varepsilon}[w]$, we apply the same procedure for $w$. \par
%	By induction over the level, we construct a sequence of cm equivalences of increasing level that starts at $X$ and ends at $P$. The concatenation of all cm equivalences of that sequence yields the desired cm equivalence $(X,f) \rightarrow (P,f')$. 
\end{proof}
\begin{Rem}
	The way we defined cm equivalence makes it a generalization of symmetry equivalence. In fact, cm equivalences are the same as symmetry equivalences, i.e. they are always trivial, if we restrict ourselves to dMfs on paths: Without loss of generality, cm equivalences of some level $a$ of a dMf on a path $(P,f)$ describe all different possibilities of how two glue two paths together with a new edge in order to obtain a path again. This means that the edge labeled $a$ can only be adjacent to the vertices that are adjacent to less than two edges, respectively, of the two old paths in $P_{a-\varepsilon}$. Thus, there are at most four possibilities for the two vertices which may be adjacent to the edge labeled $a$. All of these possibilities result in dMfs which are related to each other by reflections of the original two paths in $P_{a-\varepsilon}$. Hence, they are all symmetry-equivalent to each other.
\end{Rem}
\section{Construction of the Induced Index-Ordered DMF}\label{allatonce}
    	We address the inverse question: For any given merge tree $T$, is there a discrete Morse function $f$ on a path $P$ such that $M(P,f)\cong T$? We answer this question affirmatively by presenting an explicit construction of $P$ and two possible choices for $f$. The basic idea for the construction is to reverse-engineer the construction of the induced merge tree from \Cref{Mer}.\par
    	To start with the index-ordered case, we define two different orders on $T$. First we define a Morse order on $T$, which we call the \D{index Morse order}. Afterwards, we define the simplex order on the nodes of $T$, which we use to turn the Morse labeling induced by the index Morse order into a dMf on $P$. In \Cref{secsubcon} we will present an alternative dMf which represents $T$, namely the sublevel-connected dMf.\par
    	We will discuss in \Cref{Rel} to what extent the constructed dMf is a unique representative for $T$.

\subsection{The Index Morse Order}
To define the index Morse order, we first observe that every node $a$ of $T$ is uniquely determined by the shortest path from the root to $a$. We recall that the depth of $T$ is the maximal length of any path in $T$ that appears as the shortest path from the root to a leaf. Because $T$ is chiral, we can identify such shortest paths with certain words:
\begin{Def}\label{pw}
	Let $T$ be a merge tree of depth $n$ and let $a$ be a node of $T$. The \D{path word} corresponding to $a$ is a word $a_0a_1\dots a_n \in \{L,R,\_\}^{n+1}$ where $\_$ denotes the empty letter. If $a$ is of depth $k$, the letters $a_0\dots a_k$ are given by the chirality of the nodes belonging to the shortest path from the root to $a$. The letters $a_{k+1}\dots a_n$ are then empty.
\end{Def}
\begin{Rem}\label{maxk}
    Let $a,b$ be nodes of a merge tree $T$ and let $a_0a_1\dots a_n$ be the path word corresponding to $a$ and $b_0b_1\dots b_n$ be the path word corresponding to $b$. Then the equation $a_0=b_0=L$ always holds because we consider paths that begin at the root. Because of $a_0=b_0=L$ and because we consider finite trees, there is always a maximal $k \in \mathbb{N}$ such that $a_i=b_i$ holds for all $i\leq k$. Furthermore, the last non-empty letter of a path word is always the chirality of the considered node. 
\end{Rem}
We now define the index Morse order, which will produce an index-ordered dMf on $P$ afterwards.
\begin{Def}\label{MO}
	Let $T$ be a merge tree. 
	We define the \D{index Morse order $\leq_{io}$} on the nodes of $T$ as follows:\par
	Let $a$ and $b$ be arbitrary nodes of $T$.
	If $a$ is a leaf node and $b$ is an inner node, then we define $a\leq_{io} b$. If either both $a$ and $b$ are leaf nodes or both $a$ and $b$ are inner nodes, we consider the following:\par
Let $a_0a_1\dots a_n$ be the path word corresponding to $a$ and $b_0b_1\dots b_n$ the path word corresponding to $b$.
Furthermore, let $k\in\mathbb{N}$ be maximal such that $a_i=b_i$ for all $i\leq k$. If $a_k=b_k=\textcolor{blue}{L}/\textcolor{red}{R}$ we define $a\leq_{io} b$ if and only if one of the following cases hold:
\begin{enumerate}[a)]
	\item $\textcolor{blue}{a_{k+1}=L \text{ and } b_{k+1}=R}/\textcolor{red}{a_{k+1}=R \text{ and } b_{k+1}=L}$ 
	\item $b_{k+1}=\_$
	\item $a=b \ (\Leftrightarrow k=n)$
\end{enumerate}
\end{Def}
The index Morse order is tailor-made to induce an index-ordered dMf later on. Nonetheless, we will see in \Cref{counterex} that it is in general not the only Morse order which induces an index-ordered dMf. In \Cref{secsubcon} we will introduce a different, perhaps more natural, Morse order that is more closely related to the sublevel filtration of the induced dMf. But for now we consider an example of the index Morse order and prove that $\leq_{io}$ is actually is a Morse order.
\begin{Ex} \label{iMo}
	We consider the following merge tree $T$:
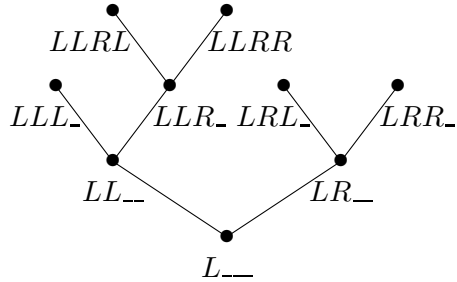
\par
	\begin{tikzpicture}
		\draw (-2.5,3.2) node {};
		\draw (5,0) node {$\bullet$};
		\draw (5,-0.4) node {\small{$L\_\_\_$}};
		\draw (5,0) -- (3.5,1);
		\draw (3.5,1) node {$\bullet$};
		\draw (3.5,0.6) node {\small{$LL\_\_$}};
		\draw (5,0) -- (6.5,1);
		\draw (6.5,1) node {$\bullet$};
		\draw (6.5,0.6) node {\small{$LR\_\_$}};
		\draw (3.5,1) -- (2.75,2);
		\draw (2.75,2) node {$\bullet$};
		\draw (2.6,1.6) node {\small{$LLL\_$}};
		\draw (3.5,1) -- (4.25,2);
		\draw (4.25,2) node {$\bullet$};
		\draw (4.5,1.6) node {\small{$LLR\_$}};
		\draw (6.5,1) -- (5.75,2);
		\draw (5.75,2) node {$\bullet$};
		\draw (5.6,1.6) node {\small{$LRL\_$}};
		\draw (6.5,1) -- (7.25,2);
		\draw (7.25,2) node {$\bullet$};
		\draw (7.5,1.6) node {\small{$LRR\_$}};
		\draw (4.25,2) -- (3.5,3);
		\draw (3.5,3) node {$\bullet$};
		\draw (3.2,2.6) node {\small{$LLRL$}};
		\draw (4.25,2) -- (5,3); 
		\draw (5,3) node {$\bullet$};
		\draw (5.3,2.6) node {\small{$LLRR$}};
		\draw (5,-0.6) node {};
	\end{tikzpicture}\\
	The path words are written underneath their corresponding nodes. The index Morse order produces the following chain of inequalities where we denote the nodes by their corresponding path words:
	\begin{align*}
		LLL\_ \lneq LLRR \lneq LLRL \lneq LRR\_ \lneq LRL\_ \lneq LLR\_ \lneq LL\_\_ \lneq LR\_\_ \lneq L\_\_\_
	\end{align*}
	The inequalities from $LLL\_$ to $LRL\_$ arise from the path words of the leaf nodes. The inequality $LRL\_ \lneq LLR\_$ holds because the node corresponding to $LRL\_$ is a leaf node and the node corresponding to $LLR\_$ is an inner node. The inequalities from $LLR\_$ to $L\_\_\_$ arise from the path words of the inner nodes.
\end{Ex}
\begin{Rem}
	By definition, the root node is always the maximal element of $(V(T),\leq_{io})$. Furthermore, the leftmost leaf node of $T$ is always the minimal element of $(V(T),\leq_{io})$.
\end{Rem}
\begin{Prop}\label{imoto}
%	The index Morse order is a total order on $V(T)$. Furthermore, it is a Morse order.
The index Morse order is a Morse order on $T$.
\end{Prop}
\begin{proof}[Sketch of Proof] The proof is a straightforward application of the definitions and involves case distinctions corresponding to the cases a), b), and c) from \Cref{MO}. 
\end{proof}
\begin{Def}\label{iML}
	We call the Morse labeling $\lambda_{io} \colon (V(T),\leq_{io}) \rightarrow \{0,1,\dots, i(T)+l(T)-1\}$ induced by the index Morse order, see \Cref{ML},  the \D{index Morse labeling} on $T$. That is, a node $c$ of $T$ is labeled with $\lambda_{io}(c)$.
\end{Def}
\begin{Ex}\label{iMlex}
	The index Morse order from \Cref{iMo} induces the following index Morse labeling:\par
		
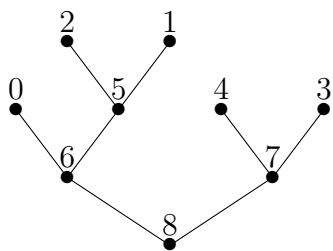
\begin{tikzpicture}[scale=0.9]
		\draw (-3.5,3.3) node {};
		\draw (5,0) node {$\bullet$};
		\draw (5,0.3) node {8};
		\draw (5,0) -- (3.5,1);
		\draw (3.5,1) node {$\bullet$};
		\draw (3.5,1.3) node {6};
		\draw (5,0) -- (6.5,1);
		\draw (6.5,1) node {$\bullet$};
		\draw (6.5,1.3) node {7};
		\draw (3.5,1) -- (2.75,2);
		\draw (2.75,2) node {$\bullet$};
		\draw (2.75,2.3) node {$0$};
		\draw (3.5,1) -- (4.25,2);
		\draw (4.25,2) node {$\bullet$};
		\draw (4.25,2.3) node {5};
		\draw (6.5,1) -- (5.75,2);
		\draw (5.75,2) node {$\bullet$};
		\draw (5.75,2.3) node {4};
		\draw (6.5,1) -- (7.25,2);
		\draw (7.25,2) node {$\bullet$};
		\draw (7.25,2.3) node {3};
		\draw (4.25,2) -- (3.5,3);
		\draw (3.5,3) node {$\bullet$};
		\draw (3.5,3.3) node {2};
		\draw (4.25,2) -- (5,3); 
		\draw (5,3) node {$\bullet$};
		\draw (5,3.3) node {1};
	\end{tikzpicture}
\end{Ex}
\subsection{The Simplex Order}
We now define the simplex order on the nodes of $T$. The simplex order will tell us which nodes of $T$ correspond to which simplices of $P$. 
	\begin{Rem}\label{yca}
	Let $T$ be a merge tree and let $a,b$ be nodes of $T$. Because $T$ is in particular a rooted binary tree, there is a unique node $p$ which is a common ancestor of $a$ and $b$ and has no descendants which are common ancestors of $a$ and $b$.    
	\end{Rem}
	\begin{Def}
	 We call the node $p$ from \Cref{yca} the \D{youngest common ancestor of $a$ and $b$}.     
	\end{Def}
\begin{Def}\label{SO}
    Let $T$ be a merge tree.
	 We define the \D{simplex order} $\preceq$ on $V(T)$ as follows: For two nodes $a$ and $b$ of $T$ we define $a\preceq b$ if and only if one of following mutually exclusive cases holds, where $p$ denotes the youngest common ancestor of $a$ and $b$:
	\begin{enumerate}
		\item $a$ is a node of the subtree with root $p_l$ and $b$ is a node of the subtree with root $p_r$.
		\item $a$ is a node of the subtree with root $b_l$ (in particular $b=p$).
		\item $b$ is a node of the subtree with root $a_r$ (in particular $a=p$).
		\item $a=b$.
		\end{enumerate}
\end{Def}
\begin{Prop}\label{sop}
	The simplex order is a total order on the nodes of $T$.
\end{Prop}
%Before we prove the proposition, we remark the following:\par
%\begin{Rem}\label{sub}
%	Let $T$ be a merge tree and let $a,b,c$ be nodes of $T$ such that $b$ is an ancestor of $c$ and a descendant of $a$. If $c$ is a node of the subtree with root \textcolor{blue}{$a_l$}/\textcolor{red}{$a_r$}, then $b$ is also a node of the subtree with root \textcolor{blue}{$a_l$}/\textcolor{red}{$a_r$}. 
%\end{Rem}
\begin{proof}[Sketch of Proof of \Cref{sop}]The proof is a bit tedious and consists of many careful case distinctions corresponding to the different cases from \Cref{SO}. Otherwise, the proof is a straightforward application of the definitions, paired with a contradiction argument here and there.
\end{proof}
We use the following definition to make the intuition of leaves being adjacent precise. This allows us to analyze the simplex order further.
\begin{Def}
     Let $T$ be a merge tree and let $a$ and $b$ be leaves of $T$. We call $a$ and $b$ \D{adjacent} if one of the following holds:
     \begin{enumerate}
         \item $a\preceq b$ and there is no leaf node $c$ of $T$ such that $a\preceq c \preceq b$ holds.
         \item $b\preceq a$ and there is no leaf node $c$ of $T$ such that $b\preceq c \preceq a$ holds.
     \end{enumerate}
\end{Def}
\begin{Lem}\label{sosc}
    Subtrees of $T$ form chains of cover relations in \mbox{$(V(T),\preceq)$}. In detail, this means the following:\par
    Let $p$ be an inner node of $T$ and let \textcolor{red}{a}/\textcolor{blue}{b} be the \textcolor{red}{leftmost}/\textcolor{blue}{rightmost} leaf of the subtree with root $p$. Then the nodes of the subtree with root $p$ in $T$ form the chain of cover relations $a\prec \dots \prec p \prec \dots \prec b$ in $(V(T),\preceq)$.
	 Moreover, any two adjacent leaves of $T$ have the property that the left one of the two adjacent leaves is covered by the youngest common ancestor of the two, whereas the right one covers the youngest common ancestor.   
\end{Lem}
\begin{proof}
    We prove the lemma inductively.	Let $p$ be a node of $T$ such that both child nodes of $p$ are leaves. It follows directly by \Cref{SO} that $p_r$ covers $p$ and $p$ covers $p_l$. That is, the subtree with root $p$ forms the chain of cover relations $a=p_l\prec p \prec p_r=b$.\par
	If $p$ is an arbitrary inner node, then by the inductive hypothesis the subtree with root \textcolor{red}{$p_l$}/\textcolor{blue}{$p_r$} forms the chain of cover relations \textcolor{red}{$a_1\prec \dots \prec p_l \prec \dots \prec b_1$}/\textcolor{blue}{$a_2 \prec \dots \prec p_r \prec \dots \prec b_2$} where \textcolor{red}{$a_1$}/\textcolor{blue}{$a_2$} is the leftmost and \textcolor{red}{$b_1$}/\textcolor{blue}{$b_2$} the rightmost leaf of the subtree with root \textcolor{red}{$p_l$}/\textcolor{blue}{$p_r$}. Since \textcolor{red}{$b_1$}/\textcolor{blue}{$a_2$} is a node of the subtree with root \textcolor{red}{$p_l$}/\textcolor{blue}{$p_r$}, it follows by case \textcolor{red}{(2)}/\textcolor{blue}{(3)} of \Cref{SO} that \textcolor{red}{$b_1 \prec p$}/\textcolor{blue}{$p \prec a_2$} holds. For all nodes $c$ of $T$ which are not nodes of the subtree with root $p$, the same case from \Cref{SO} holds for $c$ and $p$ as for $c$ and \textcolor{red}{$b_1$}/\textcolor{blue}{$a_2$}. Thus, and because \textcolor{red}{$b_1$}/\textcolor{blue}{$a_2$} is \textcolor{red}{maximal}/\textcolor{blue}{minimal} in the subtree with root \textcolor{red}{$p_l$}/\textcolor{blue}{$p_r$} by the inductive assumption, there is no node $c$ such that \textcolor{red}{$b_1\prec c \prec p$}/\textcolor{blue}{$p \prec c \prec a_2$} holds. In conclusion, \textcolor{red}{$p$ covers $b_1$}/\textcolor{blue}{$a_2$ covers $p$}. 
\end{proof}
As mentioned before, we will use the simplex order to relate the nodes of $T$ to the simplices of a path $P$. In order to do that, we now define a corresponding simplex order on the simplices of $P$.
\begin{Def}\label{so}
	Let $P$ be a path. There are two 0-simplices $p_0$ and $p_1$ in $P$ which belong only to one respective 1-simplex. For each simplex $\sigma$ of $P$ there is a unique shortest path $\gamma_\sigma$ from $p_0$ to $\sigma$. We denote the length, that is, the number of simplices, of such a path $\gamma_\sigma$ by $L(\gamma_\sigma)$. The \D{simplex order} on $P$ is defined as follows: For two simplices $\sigma$ and $\tau$ of $P$ we define $\sigma \preceq \tau$ if and only if $L(\gamma_\sigma) \leq L(\gamma_\tau)$.
\end{Def}
\begin{Lem}\label{Path}
	The simplex order on $P$ is a total order on the simplices of $P$.
\end{Lem}
\begin{proof}[Sketch of Proof]
The proof is straightforward and only uses that $P$ is a path and that the integers are linearly ordered.
\end{proof}
\begin{Rem}\label{simporcon}
    Connected subcomplexes of $P$ correspond to chains of cover relations with respect to the simplex order. Furthermore, any 1-simplex covers its left boundary 0-simplex and is covered by its right boundary 0-simplex. If one visualizes $P$ as being horizontally embedded in a plane such that $p_0$ is the leftmost 0-simplex of $P$ and $p_1$ is the rightmost 0-simplex of $P$, then for simplices $s,s'\in P$ the relation $s\prec s'$ holds if and only if $s$ is left of $s'$. This reminds us of the fact that the simplex order is only defined up to a choice of orientation.
\end{Rem}
\begin{Rem}\label{iso}
	Let $T$ be a merge tree and let $P$ be a path with $i(T)$ 1-simplices. Then we have a unique isomorphism $\phi \colon (P,\preceq) \xrightarrow{\cong} (V(T),\preceq)$ of totally ordered sets. The isomorphism $\phi$ only depends on the choice of $p_0$ and $p_1$, that is, on a choice of orientation on $P$. Choosing $p_0$ and $p_1$ the other way around would reverse the simplex order on $P$.
\end{Rem}
%\begin{proof}
%	It is immediate that $V(T)$ and $P$ have the same number of elements. 
%    Since both sets are totally ordered, finite, and have the same number of elements, there is exactly one order-preserving bijection between them. 
%	Choosing $p_0$ and $p_1$ the other way around would reverse the total order on $P$. This would result in a different isomorphism according to the inverse order. 
%\end{proof}
Before we continue with the definition of the index-ordered dMf, we consider how the simplex order can be used to classify the connected components of sublevel complexes of dMfs on paths $(P,f)$:
\begin{Prop}\label{simpconcomp}
    Let $f\colon P \rightarrow \mathbb{N}_0$ be a dMf on a path $P$.
    Then the connected components $P_c[v]$ of sublevel complexes $P_c$ of $P$ are precisely maximal sequences $\sigma\coloneqq(s_0,\dots, v, \dots, s_k)$ of simplices of $P$ such that $s_i \in P_c$ for all $i=0,\dots, k$ and $s_0\prec \dots \prec v \prec \dots \prec s_k$ is a chain of cover relations in $(P,\preceq)$.
\end{Prop}
\begin{proof}[Sketch of Proof]
	The proof is straightforward and uses \Cref{simporcon}
\end{proof}
%\begin{proof}
%    As mentioned in \Cref{simporcon}, connected subcomplexes of $P$ are given by chains $\sigma$ of cover relations in $(P,\preceq)$. Thus, it follows directly that connected components of sublevel complexes are maximal chains $\sigma$ of cover relations in $P$ that consist only of simplices of $P_c$.\par
%    For the other way around, it is clear that sequences $\sigma$ of simplices of $P_c$ form a subset of $P_c$. Because $\sigma$ is supposed to be a sequence of cover relations in $(P,\preceq)$, it is connected by \Cref{simporcon}. By the assumption that $\sigma$ is maximal, its simplices form a maximal connected subcomplex, or in other words, a connected component of $P_c$. 
%\end{proof}
\subsection{The Induced Index-Ordered DMF}
Now we explain how the simplex order can be used to construct dMfs on $P$ from Morse orders on $T$:
\begin{Def}\label{f'}
	Let $T$ be a merge tree and $P$ be a path such that the number of 1-simplices is $i(T)$ and let $\phi \colon (P,\preceq) \rightarrow (V(T),\preceq)$ be the isomorphism from \Cref{iso}.\par 
	For a Morse order $\leq$ and its induced Morse labeling $\lambda$ we define a map $f_{\lambda} \colon P\rightarrow \mathbb{N}_0$ by $f_\lambda \coloneqq \lambda \circ \phi$.
%	\begin{align*}
%		\begin{xy}
%			\xymatrix{P \ar[rr]^{f_{\lambda}} \ar[d]^{\phi} & & \mathbb{N}_0  \\
%				V(T) \ar[urr]^{\lambda} & & 
%			}
%			\end{xy}
%	\end{align*}
	The map $f_{\lambda}$ is then called the \D{dMf induced by the Morse order $\leq$} or the \D{dMf induced by the Morse labeling $\lambda$}.\par 
	In particular, the map $f_{io}\coloneqq \lambda_{io}\circ \phi$ induced by the index Morse order is called the \D{induced index-ordered dMf}.
\end{Def}
\begin{Rem}
	Although $\phi$ and $\lambda$ are order-preserving maps with respect to the previously defined total orders, the map $f_{io}$ does not respect the simplex order in general. Since $f_{io}$ is supposed to be an index-ordered dMf, it does not need to respect the simplex order. Because the map $f_{io}$ is supposed to be index-ordered, it rather needs to be compatible with face relation on $P$, which we will see to be true later on. 
\end{Rem}
\begin{Ex}
	The index Morse order from \Cref{iMo}, respectively the index Morse labeling from \Cref{iMlex}, produces the following pair $(P,f_{io})$:\par
	
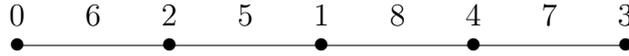
\begin{tikzpicture}
		\draw (-2.5,0.6) node {};
		\draw (1,0) -- (9,0);
		\draw (1,0.4) node {0};
		\draw (1,0) node {$\bullet$};
		\draw (3,0.4) node {2};
		\draw (3,0) node {$\bullet$};
		\draw (5,0.4) node {1};
		\draw (5,0) node {$\bullet$};
		\draw (7,0.4) node {4};
		\draw (7,0) node {$\bullet$};
		\draw (9,0.4) node {3};
		\draw (9,0) node {$\bullet$};
		\draw (2,0.4) node {6};
		\draw (4,0.4) node {5};
		\draw (6,0.4) node {8};
		\draw (8,0.4) node {7};
	\end{tikzpicture}
\end{Ex}
\begin{Prop}\label{induceddmfisdmf}
    For any given Morse order $\leq$ on any merge tree $T$ the dMf induced by $\leq$ is a dMf that has only critical cells.
\end{Prop}
\begin{proof}[Sketch of Proof]
	The proof is straightforward and uses \Cref{sosc}, \Cref{simporcon}, and property (1) of \Cref{morseorder}.
\end{proof}
%\begin{proof}
%    Because the maps $\lambda$ and $\phi$ are injective, the same is true for $f_{\lambda}$. In particular, $f_{\lambda}$ is matching and all simplices of $P$ are critical with respect to $f$.\par
%	We now prove that $f_{\lambda}$ is weakly increasing: Let $\sigma$ be a 1-simplex of $P$. By \Cref{simporcon} $\sigma$ and its boundary 0-simplices form a chain of cover relations in $(P,\preceq)$. By \Cref{sosc}, the boundary 0-simplices of $\sigma$ correspond to adjacent leaves in $T$ such that $\sigma$ is their youngest common ancestor. By property (1) of \Cref{morseorder}, it follows that $\sigma$ has a higher Morse label than its boundary 0-simplices. Thus, the map $f_{\lambda}$ is weakly increasing.\par
%    In conclusion, $f_{\lambda}$ is a dMf with only critical cells.
%\end{proof}
\begin{Rem}
    The previous proposition proves that Morse orders $\leq$ on merge trees $T$ always induce dMfs $f_{\lambda}$. It is a priori unclear though whether the induced dMf $f_{\lambda}$ induces the given merge tree $T$ as its induced merge tree $M(P,f_{\lambda})$. We prove this to be true in \Cref{MainT}. \par
    Furthermore, condition (1) from \Cref{morseorder} is necessary for $f_{\lambda}$ to be a dMf, because a violation of (1) between an inner node and a leaf would result in a violation of $f_{\lambda}$ being weakly increasing on the corresponding simplices. 
\end{Rem}
%\begin{Rem}
%	It is clear that for a given merge tree $T$ the path $P$ from this section and the path $P$ from \Cref{Con} are isomorphic as paths since they both contain $i(T)+1$ many 0-simplices and $i(T)$ many 1-simplices.
%\end{Rem}
Before we continue with the sublevel-connected dMf, we consider how the simplex order can be used to improve our understanding of sublevel complexes of dMfs on paths $(P,f_{\lambda})$ and how they are related to subtrees of $T$. The condition for this approach to be applicable is that the dMf $f_{\lambda}$ is induced by a Morse order $\leq$ on $T$ as in \Cref{f'}, which we will see to be the general case later on. We will apply this approach to the sublevel-connected case in \Cref{secsubcon} where it will be of more importance.\par
\begin{Prop}\label{concomp}
    Let $f_{\lambda}\colon P \rightarrow \mathbb{N}_0$ be a dMf on a path $P$ that is induced by a Morse order $\leq$ on $T$. Then the connected components $P_c[v]$ of sublevel complexes $P_c$ of $(P,f)$ induce subtrees of $T$ via $\phi$.
\end{Prop}
\begin{proof}[Sketch of Proof]
	It follows by \Cref{simpconcomp} and the definition of $\phi$ in \Cref{iso} that connected components of sublevel complexes $P_c[v]$ induce maximal chains of cover relations such that the corresponding simplices are of at most level $c$ in $(V(T),\preceq)$.\par
	The next step is to prove that such chains are equal to the subtree with the chain's maximum as root. The proof that the chain is contained in the subtree is straightforward. The other inclusion can be proved by contradiction, using that a node outside the subtree would contradict the property of being a chain of cover relations. 
\end{proof}
\section{The Sublevel-Connected DMF}\label{secsubcon}
As remarked in \Cref{slc} it might sometimes be more convenient to work with sublevel-connected dMfs rather than with index-ordered dMfs. In this section we introduce a slightly different version of the Morse order from \Cref{MO} to construct a sublevel-connected dMf which is shuffle-equivalent to the induced index-ordered dMf and, hence, induces the same given merge tree.  
\begin{Def}\label{scmo}
Let $T$ be a merge tree. 
We define the \D{sublevel-connected Morse order $\leq_{sc}$} on the nodes of $T$ as follows:\par
Let $a,b$ be arbitrary nodes of $T$.
Let $a_0a_1\dots a_n$ be the path word corresponding to $a$ and $b_0b_1\dots b_n$ the path word corresponding to $b$ (see \Cref{pw}).
Furthermore, let $k\in\mathbb{N}$ be maximal such that $a_i=b_i$ for all $i\leq k$. If $a_k=b_k=\textcolor{blue}{L}/\textcolor{red}{R}$ we define $a\leq_{sc} b$ if and only if one of the following cases hold:
\begin{enumerate}[a)]
	\item $\textcolor{blue}{a_{k+1}=L \text{ and } b_{k+1}=R}/\textcolor{red}{a_{k+1}=R \text{ and } b_{k+1}=L}$ 
	\item $b_{k+1}=\_$
	\item $a=b$
\end{enumerate}
\end{Def}
\begin{Rem}\label{sameorder}
	The only difference between the definition of the index Morse order and the definition of the sublevel-connected Morse order is that we do not treat leaves and inner nodes differently anymore. Thus, both orders induce the same order on inner nodes and the same order on leaves, which makes the index Morse order and the sublevel-connected Morse order shuffle-equivalent. However, the order relation between a leaf and an inner node is in general different then in the index Morse order.\par 
\end{Rem}
\begin{Rem}
    The fact that the sublevel-connected Morse order is a Morse order can be proved the same way as the corresponding statement \Cref{imoto} for the index Morse order was proved. There are just fewer case distinctions to be made for the sublevel-connected Morse order.
\end{Rem}
\begin{Lem}\label{Lemsco}
	Subtrees of $T$ form intervals in $(V(T),\leq_{sc})$.
\end{Lem}
\begin{proof}[Sketch of Proof]
	The proof is straightforward and uses the fact that all path words of a subtree $T'$ start with the same couple of letters corresponding to the root of $T'$.
\end{proof}
%\begin{proof}
%Let $T'$ be a subtree in $T$. The nodes of $T'$ correspond to path words that begin with the same couple of letters. We denote the path words of nodes $a$ of $T'$ by $a_0a_1\dots a_n$. Then there is a maximal $k \in \mathbb{N}$ such that $a_0\dots a_k$ is the same for all nodes of $T'$. In particular, the path word $a_0\dots a_k \_ \dots \_ $ corresponds to the root of $T'$.\par
%Let $b$ be a node of $T$ that is not in $T'$ and let $b_0\dots b_n$ be the corresponding path word. Then $a_i \neq b_i$ holds for at least one $i\leq k$. Thus, the order relations between $b$ and any node of $T'$ are decided by the letters $a_{i}$ and $b_{i}$ for the minimal such $i$. In particular, the relation is the same for all nodes of $T'$. Thus, all nodes $b$ which are not in $T'$ are either greater than all nodes of $T'$ or smaller than all nodes of $T'$. In conclusion, the nodes of $T'$ form a chain of cover relations because $\leq_{sc}$ is a total order.  
%\end{proof}
The induced sublevel-connected labeling $\lambda_{sc}$ and the induced sublevel-connected dMf $f_{sc}$ are defined analogously to the index-ordered case:
\begin{Def}\label{subcondmf}
	 Let $T$ be a merge tree and $P$ be a path such that the number of 1-simplices is $i(T)$.
	 We call the Morse labeling  $\lambda_{sc} \colon (V(T),\leq_{sc}) \rightarrow \{0,1,\dots, i(T)+l(T)-1\}$ induced by the sublevel-connected Morse order the \D{sublevel-connected Morse labeling} on $T$. 
	 The dMf $f_{sc}=\lambda_{sc}\circ \phi$ induced by $\lambda_{sc}$ is called the \D{induced sublevel-connected dMf}.
\end{Def}
\begin{Ex}
	Let $T$ be the merge tree from \Cref{iMo}. The sublevel-connected Morse labeling on $T$ and the induced sublevel-connected dMf are given below.\\
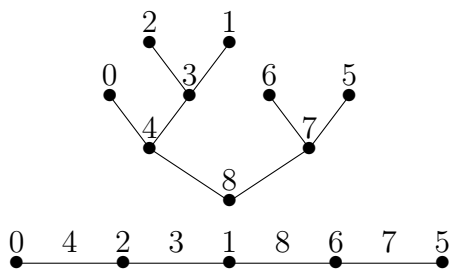
	\begin{tikzpicture}[scale=0.7]
		\draw (-5.7,3.6) node {};
		\draw (5,0) node {$\bullet$};
		\draw (5,0.4) node {8};
		\draw (5,0) -- (3.5,1);
		\draw (3.5,1) node {$\bullet$};
		\draw (3.5,1.4) node {4};
		\draw (5,0) -- (6.5,1);
		\draw (6.5,1) node {$\bullet$};
		\draw (6.5,1.4) node {7};
		\draw (3.5,1) -- (2.75,2);
		\draw (2.75,2) node {$\bullet$};
		\draw (2.75,2.4) node {$0$};
		\draw (3.5,1) -- (4.25,2);
		\draw (4.25,2) node {$\bullet$};
		\draw (4.25,2.4) node {3};
		\draw (6.5,1) -- (5.75,2);
		\draw (5.75,2) node {$\bullet$};
		\draw (5.75,2.4) node {6};
		\draw (6.5,1) -- (7.25,2);
		\draw (7.25,2) node {$\bullet$};
		\draw (7.25,2.4) node {5};
		\draw (4.25,2) -- (3.5,3);
		\draw (3.5,3) node {$\bullet$};
		\draw (3.5,3.4) node {2};
		\draw (4.25,2) -- (5,3); 
		\draw (5,3) node {$\bullet$};
		\draw (5,3.4) node {1};		
	\end{tikzpicture}\\
\begin{tikzpicture}[scale=0.7]
	\draw (-5.7,0.6) node {};
	\draw (1,0) -- (9,0);
	\draw (1,0.4) node {0};
	\draw (1,0) node {$\bullet$};
	\draw (3,0.4) node {2};
	\draw (3,0) node {$\bullet$};
	\draw (5,0.4) node {1};
	\draw (5,0) node {$\bullet$};
	\draw (7,0.4) node {6};
	\draw (7,0) node {$\bullet$};
	\draw (9,0.4) node {5};
	\draw (9,0) node {$\bullet$};
	\draw (2,0.4) node {4};
	\draw (4,0.4) node {3};
	\draw (6,0.4) node {8};
	\draw (8,0.4) node {7};
\end{tikzpicture}
\end{Ex}
There are now two things left to prove: that the map $f_{sc}$ is indeed a sublevel-connected dMf and that it induces the given merge tree $T$.
\begin{Prop}\label{scdmf}
	The induced sublevel-connected dMf $f_{sc}$ is a sublevel-connected dMf that has only critical cells.
\end{Prop}
\begin{proof}
	As a dMf induced by a Morse order, the map $f_{sc}$ is by \Cref{induceddmfisdmf} a dMf that has only critical cells.
	It is left to prove that $f_{sc}$ is sublevel-connected:\par
	By \Cref{concomp}, the connected components of sublevel complexes of $(P,f)$ induce subtrees of $T$ via $\phi$. By \Cref{Lemsco}, subtrees of $T$ form intervals in $(T,\leq_{sc})$. The sublevel-connected Morse labeling $\lambda_{sc}$ by definition maps intervals of $(T,\leq_{sc})$ to intervals of $\mathbb{N}_0$. By concatenation of these arguments, it follows that $f_{sc}=\lambda_{sc} \circ \phi$ maps connected components of sublevel complexes to intervals of $\mathbb{N}_0$, i.e.\ the map $f_{sc}$ is sublevel-connected.
\end{proof}
\begin{Theo}\label{fsciT}
	Let $T$ be a merge tree and let $P$ be a path such that the number of 1-simplices in $P$ is $i(T)$. Then $f_{io}$ and $f_{sc}$ are shuffle-equivalent where $f_{sc}$ is the induced sublevel-connected dMf and $f_{io}$ is the induced index-ordered dMf. Thus, $M(P,f_{sc})\cong M(P,f_{io})$ holds as merge trees. 
\end{Theo}
\begin{proof}
%    By \Cref{sameorder}, the restrictions of the index Morse order and the sublevel-connected Morse order to the leaves of $T$ are the same, which means that they are shuffle equivalent.
%    Thus, the index Morse labeling and the sublevel-connected Morse labeling induce the same order on the leaves of $T$. Because the induced dMfs are constructed by the concatenation of $\lambda$ with $\phi$, the induced index-ordered dMf and the induced sublevel-connected dMf induce the same order on 0-simplices. It follows analogously that both induced dMfs induce the same order on the 1-simplices of $P$. That is, the dMfs $f_{sc}$ and $f_{io}$ are shuffle-equivalent.\par
    By \Cref{sameorder}, the index Morse order and the sublevel-connected Morse order are shuffle equivalent. It follows by \Cref{shuffleequiv} and \cref{f'} the the index-ordered dMf and the sublevel-connected dMf are shuffle equivalent. Thus, $M(P,f_{sc})\cong M(P,f_{io})$ follows by \Cref{equiv}, where $f_{io}$ is the induced index-ordered dMf defined in \Cref{f'}.
\end{proof}

\section{Relationships Between Merge Trees and DMFs on Paths}\label{Rel}
In this section we want to take a look at the bigger picture again and consider how we can relate dMfs on paths and trees to merge trees. In order to do this in a structured way, we consider the different sets of dMfs and merge trees and relate them to each other by bijections that are compatible with the various notions of equivalence we introduced earlier. Afterwards, we use said bijections in order to prove that the aforementioned dMfs $f_{io}$, \Cref{f'}, and $f_{sc}$, \Cref{scdmf}, both represent the given merge tree. We also relate the sets of dMfs and merge trees mentioned here to the setting of \cite{Fiber}.  
\begin{Rem}\label{noFormanEq}
In order to obtain bijections that are compatible with the various notions of equivalence, we restrict ourselves to the case of dMfs for which all simplices are critical. Since the induced merge tree does not take matched cells into account, we would otherwise need to define a notion of equivalence similar to Forman equivalence (see \cite[Def 4.1]{JS}) of dMfs that in particular takes simple homotopy equivalences as well as combinatorial aspects of dMfs into account. This could be done by additional pre- and post-composition of the equivalences as we defined them with simple homotopy equivalences that are compatible with the given dMfs. For simplicity, we chose to leave this aspect out of this work.
\end{Rem}
We denote the set of merge trees up to isomorphism by $Mer$ and the set of dMfs on paths with only critical simplices up to symmetry equivalence by $DMF^{\text{crit}}_P$. It follows by \Cref{symequiv} that the assignment $M(\_ \ , \_)$ is well-defined on $DMF^{\text{crit}}_P$. Furthermore, the construction of the induced dMf from \Cref{f'} extends to a map which we denote by $\Phi$. Since $\phi$ is well defined up to a choice of orientation, the map $\Phi$ is in particular well defined up to symmetry equivalence. This leaves us with the following diagram:\\
\begin{figure}[H]
\centering
\begin{tikzpicture}[scale=1]
%\draw (2,4) node {};
\node (Mer) at (0,4-.3) {$Mer$};
\node (DMF) at (8,4-.3) {$DMF^{\text{crit}}_P$};
\node (MoT) at (0,2+.3) {$MoT$};
\node (MlT) at (8,2+.3) {$MlT$};
\draw[->] (0+.6,2+.1+.3) -- (8-.6,2+.1+.3) node[midway,above] {$\operatorname{iMl}$};
\draw[->] (8-.6,2-.1+.3) -- (0+.6,2-.1+.3) node[midway,below] {$\operatorname{iMo}$};
\draw[->] (DMF) -- (Mer) node[midway, above] {$M(\_ \ ,\_)$};
\draw[->] (8.1,2.3+.3) -- (8.1,4-.3-.3) node[midway,right] {$\Phi$};
\draw[->] (8-.1,4-.3-.3) -- (8-.1,2.3+.3) node[midway,left] {$M(\_ \ , \_)$};
\draw[->] (0.5,2.3+.3) -- (0.5,4-.3-.3) node[midway,right] {$\operatorname{forget}$};
\draw[->,red] (0-.3,4-.3-.3) -- (0-.3,2.3+.3) node[midway,right,red] {$\leq_{io}$};
\draw[->,red] (0-.5,4-.3-.3) -- (0-.5,2.3+.3) node[midway,left,red] {$\leq_{sc}$};
\end{tikzpicture}
\caption{Relationships Between Merge Trees and DMFs on Paths}
    \label{diag}
\end{figure}
The arrows induced by the index Morse order and the sublevel-connected Morse order are not left-inverse to the forget arrow as the following example shows. We have marked them in red because they are the only arrows that prevent the diagram from commuting completely. Nonetheless, it is obvious that the arrows induced by the two Morse orders are right-inverses of the forgetful arrow.
\begin{Ex}\label{counterex}
    Consider the following Ml trees and their corresponding Mo trees:\\
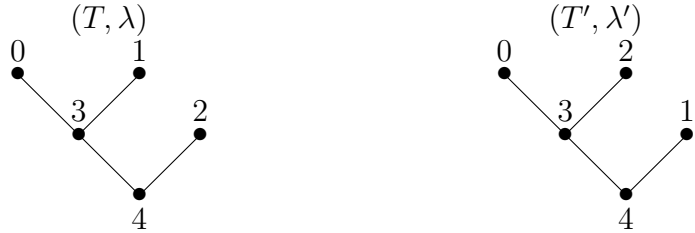
    \begin{tikzpicture}[scale=0.8]
    \draw (-4,-0.9) node {};
	\draw (1.5,-1.1) node {$(T,\lambda)$};
	\draw (2,-4) node {$\bullet$};
	\draw (2,-4.4) node {$4$};
	\draw (1,-3) node {$\bullet$};
	\draw (1,-2.6) node {$3$};
	\draw (3,-3) node {$\bullet$};
	\draw (3,-2.6) node {$2$};
	\draw (2,-2) node {$\bullet$};
	\draw (2,-1.6) node {$1$};
	\draw (0,-2) node {$\bullet$};
	\draw (0,-1.6) node {$0$};
	\draw (2,-4) -- (3,-3);
	\draw (2,-4) -- (0,-2);
	\draw (1,-3) -- (2,-2);
	\draw (1.5+8,-1.1) node {$(T',\lambda')$};
	\draw (2+8,-4) node {$\bullet$};
	\draw (2+8,-4.4) node {$4$};
	\draw (1+8,-3) node {$\bullet$};
	\draw (1+8,-2.6) node {$3$};
	\draw (3+8,-3) node {$\bullet$};
	\draw (3+8,-2.6) node {$1$};
	\draw (2+8,-2) node {$\bullet$};
	\draw (2+8,-1.6) node {$2$};
	\draw (0+8,-2) node {$\bullet$};
	\draw (0+8,-1.6) node {$0$};
	\draw (2+8,-4) -- (3+8,-3);
	\draw (2+8,-4) -- (0+8,-2);
	\draw (1+8,-3) -- (2+8,-2);
    \end{tikzpicture}\\
    It is immediate the $(T,\lambda)$ and $(T',\lambda')$ are neither isomorphic nor shuffle equivalent. Thus, $(T',\lambda')\ncong (T,\lambda)= (\operatorname{forget}(T',\lambda'),\lambda_{io}) $ holds. Nonetheless, the Ml tree $(T',\lambda')$ induces an index-ordered dMf.
\end{Ex}
But we clearly have the following.
\begin{Rem}
Let $T,T'$ be isomorphic merge trees. Then $(T,\leq_{io})\cong (T',\leq_{io})$ and $(T,\leq_{sc})\cong (T',\leq_{sc})$ holds as Mo trees. Furthermore, we have $\operatorname{forget} (T,\leq_{io})\cong T \cong \operatorname{forget}(T,\leq_{sc})$ as merge trees.
\end{Rem}

We have already seen in \Cref{moml} that the maps $\operatorname{iMo}$ and $\operatorname{iMl}$ are inverse to each other in the sense that they are bijections compatible with isomorphisms, order equivalences and shuffle-equivalences. We now show that $M(\_ \ , \_)$ and $\Phi$ are also inverse to each other up to the respective notions of equivalence.
\begin{Theo}\label{dmfml}
    The induced labeled merge tree $M(\_ \ , \_)$ and the induced dMf $\Phi$ define maps $M(\_ \ , \_)\colon DMF^{\text{crit}}_P \xleftrightarrow{} MlT \colon \Phi$ that are inverse to each other in the following sense:
    \begin{enumerate}
        \item For any dMf $(P,f)$ with only critical cells, the dMf $\Phi (M(P,f),\lambda_f)$ is symmetry-equivalent to $(P,f)$, and
        \item For any Ml tree $(T,\lambda)$, the Ml tree $M(\Phi T,f_\lambda)$ is isomorphic to $(T,\lambda)$.
    \end{enumerate}
\end{Theo}
\begin{proof}
    \begin{enumerate}
        \item Let $(P,f)$ be a dMf on a path. We construct a symmetry equivalence between $f$ and $f_{\lambda_f}$. It is given as follows: For any simplex $\sigma$ of $\Phi M(P,f)$ there is exactly one simplex $\tilde{\sigma}$ of $P$ such that $f(\tilde{\sigma})=f_{\lambda_f}(\sigma)$. This induces a bijection $\varphi\colon P \rightarrow \Phi M(P,f)$ which is compatible with $f$,$f_{\lambda_f}$ and $\operatorname{id}_\mathbb{R}$ by definition.
        But in general, the map $\varphi$ is not simplicial. This is because the simplex order on $M(P,f)$ might be different than the left/right relation on the corresponding simplices of $P$. In other words, the map $M(\_) \colon P \rightarrow M(P,f)$ is in general not compatible with the two different simplex orders. Nonetheless, the Morse labeling $\lambda_f$ induced by $f$ orders the nodes of $M(P,f)$ in the same order as their corresponding simplices of $P$. Thus, connected components of sublevel complexes of $(P,f)$ still correspond to subtrees of $M(P,f)$. Since the induced merge tree assigns the chirality of child nodes according to which connected component carries the minimal value of $f$, at each inner node of $M(P,f)$ the chirality of the two child nodes is either assigned in accordance with the left/right relation on the corresponding sublevel complex, or it is the opposite. If it is the opposite, this can be corrected by application of the reflection of the corresponding sublevel complex, that is, by application of a sublevel equivalence.
        In consequence, the difference between the right/left relation of the corresponding simplices in $P$ and the simplex order only lies in symmetry equivalences of $P$. Thus, $\varphi$ can be decomposed into a symmetry equivalence of $P$ and a simplicial isomorphism $\tilde{\varphi}$. Hence, $(\varphi,\psi)$ is a symmetry equivalence of dMfs on paths.
    \item Let $(T,\lambda)$ be an Ml tree. 
Let $c_0<c_1<\dots <c_n$ be the critical values of $f_{\lambda}$ and let $\sigma_i \in \Phi T$ such that $f_\lambda (\sigma_i)=c_i$. We recall that the induced merge tree $M$ defines in particular a bijection between the critical simplices of $\Phi T$ and the nodes of $M(\Phi T,f_\lambda)$. For any simplex $\sigma \in \Phi T$, we denote the node of $M(\Phi T,f_\lambda)$ that corresponds to $\sigma$ by $M(\sigma)$. An isomorphism $\varphi\colon(T,\lambda)\rightarrow M(\Phi T,f_\lambda)$ is given by $\varphi \coloneqq M\circ \phi^{-1}$. It is immediate that $\varphi$ is a bijection because $M$ and $\phi$ are. Furthermore, $\varphi$ is by construction compatible with the respective Morse labelings. It is only left to show that $\varphi$ is compatible with the chiral child relation and the respective roots.\par 
Consider $\sigma_n\in \Phi T$. For both trees, the simplex $\sigma_n$ corresponds to the root of the respective tree. In $M(\Phi T,f_\lambda)$ this is the case because $\sigma_n$ carries the maximal value of $f_\lambda$. In $(T,\lambda)$ this holds because $\phi(\sigma_n)$ holds the maximal Morse label $\lambda(\phi(\sigma_n))=c_n$. Thus, the map $\varphi$ maps the root of $(T,\lambda)$ to the root of $M(\Phi T,f_\lambda)$. \par
   Let $\sigma_i$ be a simplex of $\Phi T$. We now prove that $\varphi$ is compatible with the chiral child relation, that is, that \textcolor{blue}{$\varphi(\phi(\sigma_i)_l)=M(\sigma_i)_l$}/\textcolor{red}{$\varphi(\phi(\sigma_i)_r)=M(\sigma_i)_r$} holds. If $\sigma_i$ is a 0-simplex then there is nothing to show because then both $M(\sigma_i)$ and $\phi(\sigma_i)$ are leaves. Let $\sigma_i$ be a critical edge with chirality L. The case for chirality R works symmetrically to the case with chirality L. \par
   By \Cref{Mer}, the node \textcolor{blue}{$M(\sigma_i)_l$}/\textcolor{red}{$M(\sigma_i)_r$} corresponds to a critical simplex of the connected component of $\Phi T_{c_i-\varepsilon}$ that \textcolor{blue}{carries}/\textcolor{red}{does not carry} the minimal value of $f_\lambda$ on these two connected components. Furthermore, the node \textcolor{blue}{$M(\sigma_i)_l$}/\textcolor{red}{$M(\sigma_i)_r$} corresponds to the critical simplex that carries the maximal value of $f_\lambda$ of the respective connected component of $\Phi T_{c_i-\varepsilon}$. \par
   By application of \Cref{concomp}, we see that the connected components \textcolor{blue}{$\Phi T_{c_i-\varepsilon}[M^{-1}(M(\sigma_i)_l)]$}/\textcolor{red}{$\Phi T_{c_i-\varepsilon}[M^{-1}(M(\sigma_i)_r)]$} induce subtrees of $(T,\lambda)$ via $\phi$. It follows that the node \textcolor{blue}{$\phi(\sigma_i)_l$}/\textcolor{red}{$\phi(\sigma_i)_r$} is contained in the subtree that corresponds to \textcolor{blue}{$\Phi T_{c_i-\varepsilon}[M^{-1}(M(\sigma_i)_l)]$}/\textcolor{red}{$\Phi T_{c_i-\varepsilon}[M^{-1}(M(\sigma_i)_r)]$} via $\phi$ because by (2) of \Cref{morseorder} the subtree with root \textcolor{blue}{$\phi(\sigma_i)_l$ does}/\textcolor{red}{$\phi(\sigma_i)_r$ does not} carry the minimal Morse label of the subtree with root $\phi(\sigma_i)$ in $(T,\lambda)$. Furthermore, the node \textcolor{blue}{$\phi(\sigma_i)_l$}/\textcolor{red}{$\phi(\sigma_i)_r$} corresponds to the simplex that carries the maximal value of $f_\lambda$ of \textcolor{blue}{$\Phi T_{c_i-\varepsilon}[M^{-1}(M(\sigma_i)_l)]$}/\textcolor{red}{$\Phi T_{c_i-\varepsilon}[M^{-1}(M(\sigma_i)_r)]$} because by (1) of \Cref{morseorder} it carries the maximal Morse label on said subtree and because $\phi$ is order-preserving. Thus, \textcolor{blue}{$\varphi(\phi(\sigma_i)_l)=M(\sigma_i)_l$}/\textcolor{red}{$\varphi(\phi(\sigma_i)_r)=M(\sigma_i)_r$} holds.
   \end{enumerate}
\end{proof}

By considering \Cref{diag} we see that the difference between taking the induced merge tree of a dMf on a path is the same as taking its induced Ml tree and forgetting the Morse labeling. Thus, constructing a dMf that represents a given merge tree $T$ is up to symmetry equivalence the same as choosing a Morse order on $T$. This leads us to:
\begin{Theo}\label{MainT}
    Let $T$ be a merge tree and let $P$ be a path such that the number of 1-simplices in $P$ is $i(T)$. Then $T\cong M(P,f_{io})\cong M(P,f_{sc})$ holds as merge trees where $f_{io}$ denotes the induced index-ordered dMf (\Cref{f'}) and $f_{sc}$ denotes the sublevel-connected dMf (\Cref{subcondmf}).
\end{Theo}
\begin{proof}
	The statement follows by \Cref{dmfml}, \Cref{moml}, and the fact that by \Cref{MoTiso}, isomorphisms of Mo trees are in particular isomorphisms of the underlying merge trees. Furthermore, $M(P,f_{io})\cong M(P,f_{sc})$ holds by \Cref{fsciT}.
\end{proof}
%\begin{proof}
%%By \Cref{commute} and \Cref{moml},
%It is immediate that $M(P,f_{io})\cong\operatorname{forget}(\operatorname{iMo}(M(P,f_{io}),\lambda_{f_{io}}))$ holds as merge trees, where on the left-hand-side $M(P,f_{io})$ denotes the induced merge tree and on the right-hand-side $M(P,f_{io})$ denotes the induced Ml tree. By \Cref{dmfml}, $(M(P,f_{io}),\lambda_{f_{io}})\cong (T,\lambda_{io})$ holds as Ml trees. This implies $\operatorname{iMo}(M(P,f_{io},\lambda_{f_{io}}))\cong \operatorname{iMo}(T,\lambda_{io})$ as Mo trees by \Cref{moml}.
%The statement follows because by \Cref{MoTiso} isomorphisms of Mo trees are in particular isomorphisms of the underlying merge trees. Furthermore, $M(P,f_{io})\cong M(P,f_{sc})$ holds by \Cref{fsciT}. 
%\end{proof}
Using the notion of component-merge equivalence, \Cref{cmequiv}, we can extend \Cref{dmfml} to a bijection between the set of dMfs with only critical simplices on trees up to cm equivalence $DMF^{\text{crit}}_X$ and the set of Ml trees up to isomorphism $MlT$:
\begin{Theo}\label{dmfmlt}
The induced labeled merge tree $M(\_ \ , \_)$ and the induced dMf $\Phi$ define maps $M(\_ \ , \_)\colon DMF^{\text{crit}}_X \xleftrightarrow{} MlT \colon \Phi$ that are inverse to each other in the sense that:
\begin{enumerate}
	\item for any dMf $(X,f)$ with only critical cells, the dMf $\Phi (M(X,f),\lambda_f)$ is cm-equivalent to $(X,f)$, and
	\item for any Ml tree $(T,\lambda)$, the Ml tree $M(\Phi T,f_\lambda)$ is isomorphic to $(T,\lambda)$.
\end{enumerate}	
\end{Theo}
\begin{proof}
	The proof for statement (2) works exactly as in the proof for \Cref{dmfml} because symmetry equivalences are in particular cm equivalences. For (1) we apply \Cref{cmequivpath} to consider a representative of the cm equivalence class of $(X,f)$ which is a dMf on a path $(P,f')$. By \Cref{cmequiviso}, the isomorphism type of the induced Ml tree does not depend on this choice. Thus, \Cref{dmfml} implies that $\Phi (M(P,f'))$ is symmetry-equivalent to $(P,f')$. Since $(P,f')$ is cm-equivalent to $(X,f)$, so is $\Phi (M(P,f'))$.
\end{proof}
\begin{Cor}
	Let $T$ be a merge tree. By \Cref{MainT} and \Cref{dmfmlt} it follows that there are dMfs on trees $(X,f)$ such that $M(X,f)\cong T$ as merge trees. 
\end{Cor}
\begin{Cor}\label{endres}
	Applying \Cref{dmfmlt} together with \Cref{moml} and \Cref{shuffleorder} yields the result that there is a bijection $DMF_X^{crit}/_{\sim}\cong MoT$ where $\sim$ denotes order equivalence.
\end{Cor}
%\section{Discussion of Results and Comparison to Literature}
We conclude this section by discussing our results and comparing them to the results of \cite{Fiber}. Using the bijections appearing in \Cref{diag} and \Cref{dmfmlt} we replaced the question of finding dMfs on paths or arbitrary trees that represent a given merge tree $T$ by finding Morse orders on $T$ instead. This argument can be used to replace the question of classifying merge equivalence classes of dMfs on trees by classifying Morse orders. \Cref{counterex} tells us that, if a merge tree $T$ has at least three leaves, there might be different Morse orders on $T$ which are not shuffle-equivalent. That is, there are dMfs on trees $(X,f)$ contained in the image of $\Phi \circ iMl$ which induce $T$ as their merge tree but are neither isomorphic, nor symmetry-equivalent, nor cm-equivalent, nor shuffle-equivalent, nor a combination of the four to each other. We found a non-empty shuffle equivalence class of Morse orders on any merge tree, defined by either the index Morse order or the sublevel-connected Morse order, since the two have been shown to be shuffle-equivalent and, hence, merge-equivalent in \Cref{fsciT}. This allowed us to answer the aforementioned question of Johnson--Scoville affirmatively. That is, any merge tree is indeed represented by a dMf on a path, which is induced by a Morse order. Using a non-trivial cm equivalence, one can also find a dMf on a tree other than a path as a representative.\par
Furthermore, our results from \Cref{subsecmerge}, \Cref{slc}, and \Cref{Rel} allow us to structure the study of the set of merge equivalence classes of dMfs on trees using the following four notions of merge-invariant equivalences between dMfs on trees: Forman equivalence, symmetry equivalence, cm equivalence, and shuffle equivalence. \par
%\begin{itemize}
%    \item Forman equivalence,
%    \item symmetry equivalence,
%    \item cm equivalence, and
%    \item shuffle equivalence.
%\end{itemize}
We discussed in \Cref{noFormanEq} how Forman equivalences could be considered together with the other notions of equivalence and why we left Forman equivalences out of this work. In \Cref{dmfml} it becomes quite clear that passing to the induced Ml tree identifies symmetry-equivalent dMfs with each other up to isomorphism. Passing from the induced Ml tree to the induced Mo tree then identifies dMfs up to order equivalence with each other up to isomorphism of the underlying Mo tree. Last of all, passing from the induced Mo tree to the induced merge tree in particular identifies shuffle-equivalent dMfs with each other. This allows us to study the different equivalence classes separately for dMfs on paths. The more liberal notion of cm equivalence, \Cref{cmequiv}, allowed us to generalize \Cref{dmfml} to dMfs on arbitrary trees as seen in \Cref{dmfmlt}.\par
In \cite{Fiber}, the author establishes a bijection between graph-equivalence (\cite[Def 6.1]{Fiber}) classes of Morse-like (\cite[Def 6.9]{Fiber}) continuous functions on the interval that attain minima at the boundary on the one hand and isomorphism classes of chiral merge trees (\cite[Def 5.3]{Fiber}) on the other hand. \par
At first glance it might seem likely that Ml trees in the sense of \Cref{ML} and chiral merge trees in the sense of \cite[Def 5.3]{Fiber} are directly related by geometric realization and considering the corresponding abstract simplicial complex but, as mentioned in the introduction, there is a subtle difference in the construction of the induced merge tree. To be precise, the two constructions only differ in the induced chirality. In \cite[Sec 5]{Fiber} the chirality is given by which of the two merging components is the left or right one with respect to the chosen orientation on the interval. In \cite{JS} the chirality is given such that drawing the induced merge tree is compatible with the elder rule: The component with the minimal value gets the same chirality as the merged component. This means that following the same chirality leads to the oldest component.\par
This convention leads to the necessity to assume property (2) of \Cref{morseorder}: a certain compatibility between the Morse order and the chirality. The compatibility between Morse orders and the chirality implies, as seen in \Cref{symequiv}, that the induced Ml tree does not distinguish between symmetry-equivalent dMfs. If one defines induced Ml trees analogously to \cite{Fiber}, that is by inducing the chirality by a chosen orientation of the path, this notion of induced Ml trees would distinguish symmetry-equivalent dMfs on paths. Moreover, Ml trees induced by symmetry-equivalent dMfs would be related by sequences of reflections of subtrees. The definition could be as follows:
\begin{Def}\label{CML}
Let $T$ be a merge tree. A \D{Curry Morse order} is a total order $\leq$ on the nodes of $T$ such that the maximal node of any subtree is the root of said subtree. A \D{Curry Morse labeling} on a merge tree $T$ is a labeling $\lambda$ on the nodes of $T$ that induces a Curry Morse order on $T$. \par
A pair $(T,\leq)$ of a merge tree with a Curry Morse order on it is called a \D{Curry Morse ordered merge tree (CMo tree)}. A pair $(T,\lambda)$ of a merge tree with a Curry Morse labeling on it is called a \D{Curry Morse labeled merge tree (CMl tree)}.\par
Let $f\colon P \rightarrow \mathbb{R}$ be a dMf where $P$ is an oriented path. Then the CMl tree induced by $f$ is constructed as is \Cref{Mer} with the difference that the chirality is assigned according to the position of the corresponding connected components with respect to the orientation instead of according to critical values.\par
%Let $c_0<c_1<\dots < c_m$ be the critical values of $f$ that are assigned to 1-simplices. The induced CMl tree $M_C(P,f)$ is constructed by induction over these critical values in descending order. Furthermore, we label the nodes $n$ of $M_C(P,f)$ with $\lambda(n)$.\par
%		For the start of the induction we begin by creating a node $n_m$ which corresponds to the critical 1-simplex in $P$ labeled $c_m$ and setting its label $\lambda(n_m)\coloneqq c_m$ and its chirality to $L$. \par
%		For the inductive step, let $n_i$ be a node of $M(X,f)$ that corresponds to a critical 1-simplex between two 0-simplices $v$ and $w$. Define $\lambda_v\coloneqq$ $\max\{f(\sigma) \rvert \sigma \in P_{c_i-\varepsilon}[v], \sigma \text{ critical}\}$ and $\lambda_w\coloneqq\max\{f(\sigma) \vert \sigma \in P_{c_i-\varepsilon}[w], \sigma \text{ critical}\}$.
%		Two child nodes of $n_i$ are created, named $n_{\lambda_v}$ and $n_{\lambda_w}$. 
%		Then label the new nodes $\lambda(n_{\lambda_v})\coloneqq\lambda_v$ and $\lambda(n_{\lambda_w})\coloneqq\lambda_w$. If $P_{c_i-\varepsilon}[v]$ is left of $P_{c_i-\varepsilon}[w]$ with respect to the orientation on $P$, we assign $n_{\lambda_v}$ the chirality L and give $n_{\lambda_w}$ the opposite chirality. Continue the induction over the rest of the critical 1-simplices.
\end{Def}
\begin{Rem}
	CMl trees as defined above are basically the same concept as generic merge trees as defined in \cite[Def 2.2]{TrBCaba}. We stick to the name CMl trees instead of generic merge trees because we already use the term merge tree in a different way.
\end{Rem}
\begin{Ex}\label{eins}
	We consider two of the symmetry-equivalent dMfs from \Cref{symmequivalentdmf} and see how the induced CMl tree distinguishes them whereas the induced Ml trees identifies them as one:\\
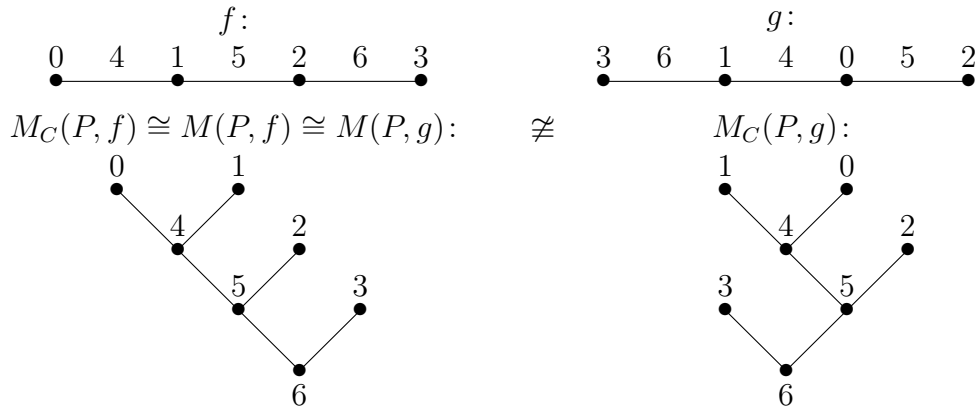
	\begin{tikzpicture}[scale=0.8]
	\draw (-2,10.2) node {};
	\draw (3,10) node {$f\colon$};
	\draw (0,9) node {$\bullet$};
	\draw (2,9) node {$\bullet$};
	\draw (4,9) node {$\bullet$};
	\draw (6,9) node {$\bullet$};
	\draw (0,9) -- (6,9);
	\draw (0,9.4) node {$0$};
	\draw (1,9.4) node {$4$};
	\draw (2,9.4) node {$1$};
	\draw (3,9.4) node {$5$};
	\draw (4,9.4) node {$2$};
	\draw (5,9.4) node {$6$};
	\draw (6,9.4) node {$3$};
	\draw (12,10) node {$g\colon$};
	\draw (9,9) node {$\bullet$};
	\draw (11,9) node {$\bullet$};
	\draw (13,9) node {$\bullet$};
	\draw (15,9) node {$\bullet$};
	\draw (9,9) -- (15,9);
	\draw (9,9.4) node {$3$};
	\draw (10,9.4) node {$6$};
	\draw (11,9.4) node {$1$};
	\draw (12,9.4) node {$4$};
	\draw (13,9.4) node {$0$};
	\draw (14,9.4) node {$5$};
	\draw (15,9.4) node {$2$};
\end{tikzpicture}\\
\begin{tikzpicture}[scale=0.8]
	\draw (-2,8) node {};
	\draw (3,8) node {$M_C(P,f) \cong M(P,f)\cong M(P,g)\colon$};
	\draw (8,8) node {$\ncong$};
	\draw (12,8) node {$M_C(P,g)\colon$};
	\draw (1,7) node {$\bullet$};
	\draw (3,7) node {$\bullet$};
	\draw (2,6) node {$\bullet$};
	\draw (4,6) node {$\bullet$};
	\draw (3,5) node {$\bullet$};
	\draw (5,5) node {$\bullet$};
	\draw (4,4) node {$\bullet$};
	\draw (1,7+.4) node {$0$};
	\draw (3,7+.4) node {$1$};
	\draw (2,6+.4) node {$4$};
	\draw (4,6+.4) node {$2$};
	\draw (3,5+.4) node {$5$};
	\draw (5,5+.4) node {$3$};
	\draw (4,4-.4) node {$6$};
	\draw (4,4) -- (1,7);
	\draw (4,4) -- (5,5);
	\draw (3,5) -- (4,6);
	\draw (2,6) -- (3,7);
	\draw (11,7) node {$\bullet$};
	\draw (13,7) node {$\bullet$};
	\draw (12,6) node {$\bullet$};
	\draw (14,6) node {$\bullet$};
	\draw (11,5) node {$\bullet$};
	\draw (13,5) node {$\bullet$};
	\draw (12,4) node {$\bullet$};
	\draw (11,7+.4) node {$1$};
	\draw (13,7+.4) node {$0$};
	\draw (12,6+.4) node {$4$};
	\draw (14,6+.4) node {$2$};
	\draw (11,5+.4) node {$3$};
	\draw (13,5+.4) node {$5$};
	\draw (12,4-.4) node {$6$};
	\draw (12,4) -- (11,5);
	\draw (12,4) -- (14,6);
	\draw (13,5) -- (11,7);
	\draw (12,6) -- (13,7);
	\end{tikzpicture}
\end{Ex}
The notion of CMl trees is related to the notion of chiral merge trees in the sense of \cite{Fiber} by the interplay between abstract and geometrical simplicial complexes. In detail, the bijection is given as follows:
\begin{Con}
Let $(T,\lambda)$ be a CMl tree. We define a chiral merge tree $\lvert(T,\lambda)\rvert$ associated to $(T,\lambda)$ as follows:\par
The compact rooted tree is given by the geometric realization $\lvert T \rvert$. We attach a distinguished edge $e_\infty$ to the vertex which corresponds to the root of $T$ in order to obtain a cell complex which we will by abuse of notation also refer to as $\lvert T \rvert$. The map $\pi\colon \lvert T \rvert \rightarrow \mathbb{R}$ is given by $\lambda$ on vertices, and by a linear extension of $\lambda$ on edges.\par 
%To make the linear extensions on edges between parent nodes $p$ and child nodes $c$ unique, we assume that every edge in $\lvert T \rvert$ has length one and choose extensions of the form $\lambda(c) + t\cdot (\lambda (p)-\lambda(c))$, where $t\in [0,1]$. The chirality on the edges in $\lvert T \rvert$ is given by the chirality of the child nodes adjacent to the corresponding edges in $T$.\par
For the other way around let $\pi\colon T \rightarrow \mathbb{R}$ be a chiral merge tree in the sense of \cite[Def 5.3]{Fiber}. We define an Ml tree $\operatorname{abs}(T)$ associated to $T$ as follows:\par
We take the 0-skeleton $T_0$ as the vertex set and the 1-skeleton $(T\setminus \{e_\infty\})_1$ as the set of edges. We define the node which corresponds to $v_\infty$ to be the root of $\operatorname{abs}(T)$. The labeling $\lambda$ is given by $\pi$. \par 
%The chiral child relation is induced by the chirality of the edges between the corresponding parent and child nodes in $T$. The labeling $\lambda$ is directly induced by $\pi$. 
The proof that the two constructions are inverse to each other is straightforward.
\end{Con}
Furthermore, there is a similar bijection between the two notions of Morse functions:
\begin{Con}
Let $(P,f)$ be a dMf with only critical cells on an oriented path. We define a Morse-like function $\tilde{f}$ on the interval which attains minima at the boundary as follows:\par
Let $k+1$ be the number of 0-simplices in $P$. We denote the 0-simplices of $P$ by $n_0,n_1,\dots , n_k$ from left to right with respect to the given orientation. Then we define $\tilde{f}(\frac{i}{k})\coloneqq f(n_i)$ for $i\in \{0,1,\dots,k\}$. We denote the 1-simplices of $P$ by $e_1, \dots, e_k$, again according to the given orientation. Then we define $\tilde{f}(\frac{2i-1}{2k})\coloneqq f(e_i)$ for $i\in \{0,1,\dots,k\}$. We define $\tilde{f}$ on the rest of the interval as the linear extension. This makes $\tilde{f}$ a distinct-valued PL function that attains minima at the boundary. Hence, $\tilde{f}$ is Morse-like.\par
For the other way around let $f\colon I \rightarrow \mathbb{R}$ be a Morse-like function that attains minima at the boundary. We define a path $P$ and a dMF $f'$ on $P$ as follows:\par
Let $n_0,n_1,\dots , n_k$ be the local minima of $f$ ordered by the orientation on $I$, so in particular $n_0=0$ and $n_k=1$.
We define $P$ to be a path with $k+1$ simplices of dimension 0. We choose one of the endpoints to be denoted by $s_0$ and the one by $s_k$. We denote the other 0-simplices such that their indices are in accordance with their position in the simplex order (\Cref{so}), making $s_0$ the minimal simplex with respect to the simplex order. The dMf $f'$ is defined by $f'(s_i)\coloneqq f(n_i)$ on 0-simplices. Let $c_1,\dots ,c_k$ be the local maxima of $f$, ordered in accordance to the orientation on $I$, and let $e_1,\dots, e_k$ be the 1-simplices of $P$, ordered in accordance to the aforementioned simplex order on $P$. Then $f'$ is defined by $f'(e_i)\coloneqq f(c_i)$ on 1-simplices.\par
It is easy to check that the two given constructions are inverse to each other.
\end{Con}
A theorem which, analogously to \Cref{dmfml}, defines a pair of inverse bijections $M_C(\_ \ ,  \_) \colon ^+DMF^{\text{crit}}_P\leftrightarrow CMlT\colon \Phi $ can be proved similarly to the proof of \Cref{dmfml}. The difference is that the induced CMl tree keeps track of symmetry equivalences in the sense that symmetry equivalences of dMfs induce reflections at roots of subtrees on the induced CMl trees. The map $\Phi$ can be defined the same way as before. Alternatively, one could define the bijection $\Phi$ as the composition $ \operatorname{abs} \circ \Psi^{-1}\circ\lvert \_ \rvert$ in the following diagram:
\begin{figure}[H]
	\centering
	\begin{tikzpicture}
	\draw (2,4) node {};
	\node (Mer) at (0,4) {$Mer$};
	\node (DMF) at (4,4) {$DMF^{\text{crit}}_P$};
	\node (MoT) at (0,2) {$MoT$};
	\node (MlT) at (4,2) {$MlT$};
	\draw[->] (0+.6,2+.1) -- (4-.6,2+.1) node[midway,above] {$\operatorname{iMl}$};
	\draw[->] (4-.6,2-.1) -- (0+.6,2-.1) node[midway,below] {$\operatorname{iMo}$};
	\draw[->] (DMF) -- (Mer) node[midway, above] {$M(\_ \ ,\_)$};
	\draw[->] (4.1,2.3) -- (4.1,4-.3) node[midway,right] {$\Phi$};
	\draw[->] (4-.1,4-.3) -- (4-.1,2.3) node[midway,left] {$M(\_ \ , \_)$};
	\draw[->] (0.5,2.3) -- (0.5,4-.3) node[midway,right] {$\operatorname{forget}$};
	\draw[->,red] (0-.3,4-.3) -- (0-.3,2.3) node[midway,right,red] {$\leq_{io}$};
	\draw[->,red] (0-.5,4-.3) -- (0-.5,2.3) node[midway,left,red] {$\leq_{sc}$};
	\node (DMF+) at (8,4) {$^+DMF^{\text{crit}}_P$};
	\node (Morselike) at (12,4) {$\mathcal{M}$};
	\node (CMlT) at (8,2) {$CMlT$};
	\node (chiral) at (12,2) {$\mathcal{X}$};
\draw[->] (0+4+.9,4+.1) -- (4+4-1,4+.1) node[midway,above] {$o_{JS}$};
\draw[->] (4+4-1,4-.1) -- (0+4+.9,4-.1) node[midway,below] {$/_{symm}$};
	\draw[->] (8.1,2.3) -- (8.1,4-.3) node[midway,right] {$\Phi$};
	\draw[->] (8-.1,4-.3) -- (8-.1,2.3) node[midway,left] {$M_C(\_ \ , \_)$};
	\draw[->] (0+8+.8,2+.1) -- (4+8-.4,2+.1) node[midway,above] {$\lvert \_ \rvert$};
	\draw[->] (4+8-.4,2-.1) -- (0+8+.8,2-.1) node[midway,below] {$\operatorname{abs}$};
	\draw[->] (0+8+1,4+.1) -- (4+8-.4,4+.1) node[midway,above] {$\lvert \_ \rvert$};
	\draw[->] (4+8-.4,4-.1) -- (0+8+1,4-.1) node[midway,below] {$\operatorname{abs}$};
	\draw[->] (12.1,2.3) -- (12.1,4-.3) node[midway,right] {$\Psi^{-1}$};
	\draw[->] (12-.1,4-.3) -- (12-.1,2.3) node[midway,left] {$\Psi$};
	\draw[->] (0+4+0.6,2+.1) -- (4+4-.8,2+.1) node[midway,above] {$i$};
	\draw[->] (4+4-.8,2-.1) -- (0+4+0.6,2-.1) node[midway,below] {$JS$};
	\end{tikzpicture}
	\caption{Relationship to the continuous case}
	\label{diag2}
\end{figure}
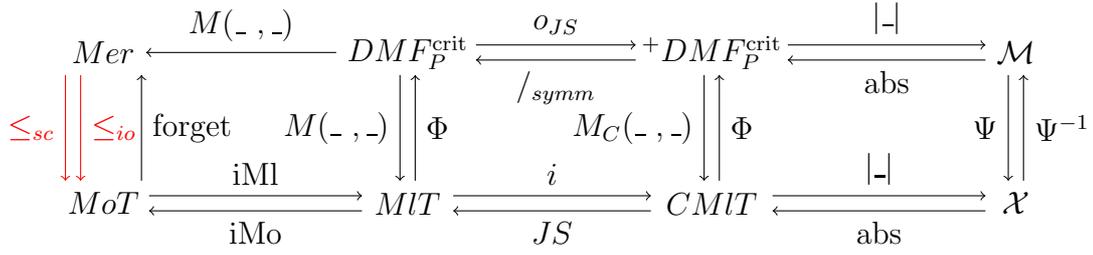
Here $\Psi$ denotes the bijection between the set of Morse-like functions on the interval $\mathcal{M}$ and the set of chiral merge trees $\mathcal{X}$ from \cite[Cor 6.11]{Fiber}. The map $i$ is the inclusion induced by considering Ml trees as CMl trees. Ml trees can be considered as a special kind of CMl trees because they only differ in the additional property (2) of \Cref{morseorder} which does not need to hold for CMl trees. The map $/_{symm}$ is the quotient map that identifies symmetry equivalent dMFs. The map $o_{JS}$ is defined as the composition $\Phi \circ i \circ M(\_ \ , \_)$. This definition coincides with choosing the representative of a dMf $(P,f)$ with respect to symmetry equivalence such that the simplex order on $P$ is compatible with $f$ in the following way: At each critical edge $e$ with $f(e)=c$ the connected component of $P_{c-\varepsilon}$ that corresponds to \textcolor{blue}{$M(e)_l$}/\textcolor{red}{$M(e)_r$} is \textcolor{blue}{left}/\textcolor{red}{right} of the edge $e$ with respect to the orientation induced by the simplex order. Here we recall that $M(e)$ denotes the inner node of $M(P,f)$ that corresponds to the critical edge $e$ and that \textcolor{blue}{$M(e)_l$}/\textcolor{red}{$M(e)_r$} denotes the  \textcolor{blue}{left}/\textcolor{red}{right} child node of $M(e)$. In \Cref{eins} we have $o_{JS}(P,g)=(P,f)$ where $(P,f)$ is oriented from left to right.\par
 The map $JS$ is defined as $JS \coloneqq M(\_ \ ,\_) \circ /_{symm} \circ \Phi$. This definition coincides with division by the equivalence relation generated by reflections of subtrees.
 \begin{Ex}
 	We consider how a sequence of reflections of subtrees maps the CMl trees from \Cref{eins} to each other. Here, we denote by $a_\lambda$ the reflection of the subtree with root labeled $\lambda$.\\
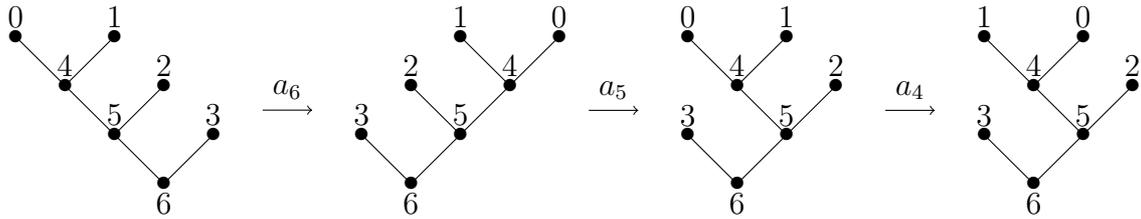
 	\begin{tikzpicture}[scale=0.65]
 	\draw (2,7.6) node {};
 	\draw (2,7) node {$\bullet$};
 	\draw (4,7) node {$\bullet$};
 	\draw (3,6) node {$\bullet$};
 	\draw (5,6) node {$\bullet$};
 	\draw (4,5) node {$\bullet$};
 	\draw (6,5) node {$\bullet$};
 	\draw (5,4) node {$\bullet$};
 	\draw (2,7+.4) node {$0$};
 	\draw (4,7+.4) node {$1$};
 	\draw (3,6+.4) node {$4$};
 	\draw (5,6+.4) node {$2$};
 	\draw (4,5+.4) node {$5$};
 	\draw (6,5+.4) node {$3$};
 	\draw (5,4-.4) node {$6$};
 	\draw (5,4) -- (2,7);
 	\draw (5,4) -- (6,5);
 	\draw (4,5) -- (5,6);
 	\draw (3,6) -- (4,7);
 	\draw[->] (7,5.5) -- (8,5.5) node[midway,above] {$a_6$};
 	\draw (11,7) node {$\bullet$};
 	\draw (13,7) node {$\bullet$};
 	\draw (10,6) node {$\bullet$};
 	\draw (12,6) node {$\bullet$};
 	\draw (9,5) node {$\bullet$};
 	\draw (11,5) node {$\bullet$};
 	\draw (10,4) node {$\bullet$};
	\draw (11,7+.4) node {$1$};
	\draw (13,7+.4) node {$0$};
	\draw (10,6+.4) node {$2$};
	\draw (12,6+.4) node {$4$};
	\draw (9,5+.4) node {$3$};
	\draw (11,5+.4) node {$5$};
	\draw (10,4-.4) node {$6$};
	\draw (10,4) -- (9,5);
	\draw (10,4) -- (13,7);
	\draw (11,5) -- (10,6);
	\draw (12,6) -- (11,7);
	\end{tikzpicture}
	\begin{tikzpicture}[scale=0.65]
 	\draw[->] (0,0.5) -- (1,0.5) node[midway,above] {$a_5$};
 	\draw (4,7-5) node {$\bullet$};
 	\draw (2,7-5) node {$\bullet$};
 	\draw (5,6-5) node {$\bullet$};
 	\draw (3,6-5) node {$\bullet$};
 	\draw (2,5-5) node {$\bullet$};
 	\draw (4,5-5) node {$\bullet$};
 	\draw (3,4-5) node {$\bullet$};
 	\draw (4,7+.4-5) node {$1$};
 	\draw (2,7+.4-5) node {$0$};
 	\draw (5,6+.4-5) node {$2$};
 	\draw (3,6+.4-5) node {$4$};
 	\draw (2,5+.4-5) node {$3$};
 	\draw (4,5+.4-5) node {$5$};
 	\draw (3,4-.4-5) node {$6$};
 	\draw (3,4-5) -- (2,5-5);
 	\draw (3,4-5) -- (5,6-5);
 	\draw (4,5-5) -- (2,7-5);
 	\draw (3,6-5) -- (4,7-5);
 	\draw[->] (6,.5) -- (7,.5) node[midway,above] {$a_4$};
 	\draw (8,7-5) node {$\bullet$};
 	\draw (10,7-5) node {$\bullet$};
 	\draw (9,6-5) node {$\bullet$};
 	\draw (11,6-5) node {$\bullet$};
 	\draw (8,5-5) node {$\bullet$};
 	\draw (10,5-5) node {$\bullet$};
 	\draw (9,4-5) node {$\bullet$};
 	\draw (8,7+.4-5) node {$1$};
 	\draw (10,7+.4-5) node {$0$};
 	\draw (9,6+.4-5) node {$4$};
 	\draw (11,6+.4-5) node {$2$};
 	\draw (8,5+.4-5) node {$3$};
 	\draw (10,5+.4-5) node {$5$};
 	\draw (9,4-.4-5) node {$6$};
 	\draw (9,4-5) -- (8,5-5);
 	\draw (9,4-5) -- (11,6-5);
 	\draw (10,5-5) -- (8,7-5);
 	\draw (9,6-5) -- (10,7-5);
 	\end{tikzpicture}
 \end{Ex}
The proof that the two possible definitions for $\Phi$ coincide and that the diagram from \Cref{diag2} commutes are straightforward. \par
 As a result, the map $\Phi$ can be seen as a discrete version of the map $\Psi$ from \cite[Cor 6.11]{Fiber}. Moreover, moving from the setting of \cite{Fiber} to the setting of \cite{JS} basically means to divide by symmetry equivalences which allows the authors of \cite{JS} to generalize the construction of the induced merge tree to dMfs on trees.
 In order to extend our generalization \Cref{dmfmlt} to the oriented case, one would need to find a notion of orientation-preserving cm equivalences which distinguishes symmetry-equivalent dMfs. 
%Furthermore, we conjecture that the correspondence between dMfs, Ml trees and Mo trees can be categorified. That is, there are categories of dMfs, Ml trees and Mo trees that have the aforementioned types of equivalences as isomorphisms and the aforementioned bijections define equivalences of categories.
\section{Further Directions and Possible Applications}
In this section we want to take a look at possible applications of our results. The structured overview of the different notions of equivalence of discrete Morse functions and their connections to each other might be useful to explore the space of discrete Morse functions on a given simplicial complex. Even though the construction of the induced merge tree from \cite{JS} does not easily extend to arbitrary simplicial complexes, the notions of equivalence from this article do. Hence, one could try to use e.g. the notion of symmetry equivalence to structure the space of discrete Morse functions in a nice way, i.e. into orbits of a groupoid action. In a second step, one could then try to define the space of merge trees as a quotient of the space of discrete Morse functions on some ``large enough'' simplicial complex. With such a construction, one could assemble the results from this article and the results from \cite{TrBCaba} into an analysis of a larger instance of the persistence map. \par
Aside from this, one could try to generalize and enhance the construction of the induced merge tree from \cite{JS} to arbitrary simplicial complexes. Then, in a second step, one could try to use the induced merge tree to find possible cancellations of pairs of critical simplices. This way, the induced merge tree might be helpful to optimize discrete Morse functions.\par
Furthermore, one could study the possible Morse orders on a given merge tree. With a classification of all Morse orders on a given merge tree $T$, one could classify all discrete Morse functions on any given tree that induce $T$ as the induced merge tree with the help of \Cref{f'} and cm equivalences. \par
Ultimately, the proofs of \Cref{dmfmlt} and all the lemmas that lead to it describe exactly which information is lost when considering the induced Ml tree instead of the original dMf. This knowledge might be useful for applications in TDA because it tells the user which kind of features will not be seen by the induced merge tree, and hence, by the persistent zeroth homology and the barcode. Moreover, the knowledge about the exact lost information can be used to enhance the induced merge tree with extra structure such that it no longer disregards certain desired information.
\bibliographystyle{alpha}
\bibliography{References}
\end{document}